\documentclass[]{siamart}
\pdfoutput=1
\newcommand{\TheTitle}{Multiscale Finite Element Modeling of Nonlinear Magnetoquasistatic 
Problems using Magnetic Induction Conforming Formulations}
\newcommand{\TheTitled}{Multiscale FE Modeling of Magnetoquasistatic Problems}
\newcommand{\TheAuthors}{I. Niyonzima, R.~V. Sabariego, P. Dular, K. Jacques and C. Geuzaine}

\headers{\TheTitled}{\TheAuthors}

\title{{\TheTitle}}

\author{
    I. Niyonzima 
    \and
    R.~V. Sabariego 
    \and
    P. Dular 
    \and
    K. Jacques
    \and
    C. Geuzaine
}

\usepackage[utf8]{inputenc}
\usepackage{amsfonts}
\usepackage{amssymb}
\usepackage{framed}
\usepackage{morefloats}
\usepackage{anyfontsize}
\usepackage{import}
\usepackage{graphicx}
\usepackage{amsmath}
\usepackage{algorithmicx}
\usepackage{algorithm}
\usepackage{algpseudocode}
\usepackage{pdfcomment}
\usepackage[final]{changes}

\usepackage{tikz}

\usetikzlibrary{matrix}
\usepackage{tcolorbox}
\usepackage{collectbox}

\algnewcommand\algorithmicinput{\textbf{INPUT:}}
\algnewcommand\INPUT{\item[\algorithmicinput]}
\algnewcommand\algorithmicoutput{\textbf{OUTPUT:}}
\algnewcommand\OUTPUT{\item[\algorithmicoutput]}

\definechangesauthor[name={Ruth V. Sabariego}, color=green]{rvs}
\setremarkmarkup{\footnote{#1:\textcolor{Changes@Color#1}{#2}}}

\renewcommand{\vec}[1]{\boldsymbol{#1}}

\newcommand{\bx}{\mbox{\boldmath$x$}}
\newcommand{\by}{\mbox{\boldmath$y$}}
\newcommand{\bn}{\mbox{\boldmath$n$}}
\newcommand{\bu}{\mbox{\boldmath$u$}}
\newcommand{\bv}{\mbox{\boldmath$v$}}

\newcommand{\ba}{\mbox{\boldmath$a$}}
\newcommand{\be}{\mbox{\boldmath$e$}}
\newcommand{\bh}{\mbox{\boldmath$h$}}

\newcommand{\bb}{\mbox{\boldmath$b$}}
\newcommand{\bj}{\mbox{\boldmath$j$}}

\newcommand{\bff}{\mbox{\boldmath$f$}}
\newcommand{\bM}{\mbox{\boldmath$M$}}
\newcommand{\bW}{\mbox{\boldmath$W$}}
\newcommand{\bE}{\boldsymbol{\mathcal{E}}}
\newcommand{\bH}{\boldsymbol{\mathcal{H}}}

\newcommand{\bB}{\boldsymbol{\mathcal{B}}}

\newcommand{\bR}{\boldsymbol{\mathcal{R}}}

\newcommand{\bF}{\boldsymbol{\mathcal{F}}}

\newcommand{\ah}{\mbox{$\alpha$}}
\newcommand{\bah}{\mbox{\boldmath$\alpha$}}

\newcommand{\Grad}[2][]  {\text{\GradSymb}{#1}\,{#2}}
\newcommand{\Curl}[2][]  {\text{\CurlSymb}{#1}\,{#2}}
\newcommand{\Div}[2][]  {\text{\DivSymb}{#1}\,{#2}}

\numberwithin{equation}{section}

\numberwithin{thm}{section}

\newtheorem{prob}{Problem}
\numberwithin{prob}{section}
\newtheorem{ass}{Assumption}

\newcommand{\pvec}[2]  {{#1}\times{#2}}

\newcommand{\GradSymb}{\text{\bf grad}}
\newcommand{\CurlSymb}{\text{\bf curl}}
\newcommand{\DivSymb}{\text{div}}
\newcommand{\isur}[3][]{\langle#2,#3\rangle#1}
\newcommand{\ivol}[3][]{(#2,#3)#1}

\newcommand{\Ltwo}[2][]    {L^2#1(#2)}

\newcommand{\LLtwo}[2][]   {\boldsymbol{L}^2#1(#2)}

\newcommand{\Hone}[2][]    {H^1#1(#2)}

\newcommand{\Hcurl}[2][]   {\boldsymbol{H}#1(\CurlSymb;#2)}
\newcommand{\Hdiv}[2][]    {\boldsymbol{H}#1(\DivSymb;#2)}

\newcommand{\Hcurlzero}[2][]   {\boldsymbol{H}#1(\textbf{curl} \, \boldsymbol{0};#2)}
\newcommand{\Hdivzero}[2][]    {\boldsymbol{H}#1(\textbf{div} \, 0;#2)}

\makeatletter
\newcommand{\mybox}{%
    \collectbox{%
        \setlength{\fboxsep}{4pt}%
        \setlength{\fboxrule}{1.5pt}%
        \fbox{\BOXCONTENT}%
    }%
}
\makeatother

\let\cite=\cite

\begin{document}

\maketitle

\begin{abstract}
    In this paper we develop magnetic induction conforming
    multiscale formulations for magnetoquasistatic problems involving periodic
    materials. The formulations are derived using the
    periodic homogenization theory and applied within a heterogeneous
    multiscale approach. Therefore the fine-scale problem
    is replaced by a macroscale problem defined on a coarse mesh
    that covers the entire domain and many mesoscale
    problems defined on finely-meshed small areas around some points of
    interest of the macroscale mesh (e.g. numerical quadrature points). The
    exchange of information between these macro and
    meso problems is thoroughly explained in this paper. For the sake
    of validation, we consider \added{a} two-dimensional \replaced{geometry of}{geometries: a
    laminated iron stack and} an idealized periodic soft magnetic composite.
\end{abstract}



\section{Introduction}
\label{sec:introduction}

The use of numerical methods for solving electromagnetic problems is nowadays
widespread. Indeed, analytical solutions of Maxwell's equations are not always
available when facing the complexity of real-life devices with complicated
geometries and materials exhibiting a possibly nonlinear or hysteretic
behaviour.  In this paper we are interested in multiscale magnetoquasistatic
(MQS) problems. These problems arise from Maxwell's equations when the
wavelength of the exciting source is much greater than the size of the structure
so that the displacement currents can be neglected. This is the model that
describes the physics of most electric power systems: electric generators,
motors and transformers.

The finite element (FE) method is a frequently-used numerical method for solving
MQS problems for its easiness to handle problems involving both nonlinearities
and complex geometries. To this end, a mesh of the structure is generated and
Maxwell's equations are weakly verified on average on elements of the mesh,
which is ensured by integrating these equations elementwise.  If the problem is
well-posed, the finer the mesh, the more accurate the numerical solution.

Soft ferrites, lamination stacks and soft magnetic composites (SMC) are
multiscale materials used in MQS applications. For instance, soft ferrites help
reducing the magnetic losses in high-frequency transformers; the cores of
electrotechnical devices are laminated to limit the eddy current losses; and the
SMCs ease the manufacturing of three-dimensional paths in electrical machines.

For problems involving such multiscale materials, the application of classical
numerical methods such as the FE method becomes prohibitive in terms of the
computational resources (time and memory) storage whence the use of
homogenization and multiscale methods.  
\added{
Using these methods, the multiscale problem is replaced by the homogenized problem 
defined on the homogeneous domain with a slowly varying fields. 
}
The performance of homogenization and multiscale methods
for MQS problems can be compared by evaluating their ability to
\begin{itemize}
\item derive a homogenized problem that can be easily solved;
\item handle nonlinearities; 
\item deal with materials with complex microstructures;
\item deal with partial differential equations involving $\CurlSymb$ operators;
\item compute global quantities such as the eddy currents or magnetic losses.
\item recover local fields at critical points of interest;
\end{itemize}

The first homogenization approach used to analytically characterize properties
of composites materials was based on mixing
rules \cite{maxwell-garnett-mixing-04,Sihvola-mixing-99}. More elaborate
theoretical methods such as the asymptotic expansion method
\cite{bensoussan-ahm-11}, the
G-convergence \cite{murat-gconv-77,tartar-homog-09}, the $\Gamma$-convergence 
\cite{DeGiorgio-gammaconv-84,braides-gamma-02,DalMaso-gammaconv-93}, the two-scale
convergence \cite{nguetseng-tsh-89,visintin-tsh-08} and the periodic unfolding
methods \cite{cioranescu-pum-02,cioranescu-pum-08} allow to construct the
homogenized problem and determine the associated constitutive laws. Equations
resulting from these methods can be used to develop multiscale methods. A
non-exhaustive list of these multiscale methods include the mean-field
homogenization method \cite{laurence-mfh-10,corcolle-thesis-09},
the multiscale finite element method--MsFEM
\cite{hou-msfem-97,efendiev-msfem-04-1}, the variational multiscale method--VMS
\cite{brezzi-vms-97,juanes-vms-05} and the heterogeneous multiscale
method--HMM \cite{e-hmm-07,Abdulle:2009p1255,e-hmm-11}.  In electromagnetism
such methods have been developed mainly for materials with
linear \cite{Bossavit1994,Bossavit1996,gyselinck-homogenization-04,meunier-chm-10,bottauscio-msfem-13,bottauscio-msfem-13-1}
and
nonlinear~\cite{gyselinck-homogenization-06,belkadi-homogenization-09,bottauscio-chm-08}
magnetic material laws. While some preliminary results concerning
electromagnetic hysteresis can be found in
\cite{12.cefc.sabariego.homogLamHyst}, there is to date no generic multiscale
method able to accurately handle hysteretic materials in complex geometrical
configurations.

In this paper we develop such a multiscale method to treat magnetoquasistatic
problems involving multiscale materials that can exhibit linear, nonlinear or
hysteretic behaviour
\added{
with the main focus on the development of weak formulations for the homogenized 
problem. Using results from the theory of homogenization for nonlinear electromagnetic 
multiscale problem obtained by Visintin, we develop the magnetic vector potential 
formulations for the multiscale, the macroscale and the mesoscale problems. 
The formulations are then validated on simple 2D geometry.
}
The multiscale method is inspired by the HMM and is based on the
scale separation assumption $\varepsilon \ll 1$ where $\varepsilon = l/L$ is the
ratio between the smallest scale $l$ and the scale of the material or the
characteristic length of external loadings $L$. The fine-scale problem is
replaced by a macroscale problem defined on a coarse mesh covering the entire
domain and many mesoscale problems that are defined on small, finely meshed
areas around some points of interest of the macroscale mesh (e.g. numerical
quadrature points). The transfer of information between these problems is
performed during the \textit{upscaling} and the \textit{downscaling} stages that
will be detailed hereafter. 

The paper comprises five sections. 
In Section \ref{section:derivation_mqs} we derive the MQS multiscale and homogenized
problems from the multiscale problem that was studied by Visintin
in \cite{visintin-tsh-06-b, visintin-tsh-08}.
In Section \ref{section:mqs} we derive the weak forms of the multiscale MQS problem.
Section \ref{section:models} deals with the multiscale weak formulations for homogenized 
MQS problems. Starting from the distributional equations that govern the fields of the 
MQS homogenized problem we develop magnetic vector potential formulations for the 
macroscale and the mesoscale problems. 
Scale transitions are also thoroughly investigated. 
Section \ref{section:tests} concerns the application of the theory to a simple but 
representative two-dimensional problem: the modeling of a soft magnetic composite. 
Conclusions are drawn in the last section.
\vspace{10mm}
{\color{black}
\section{Derivation of the homogenized magnetoquasistatic problem}
\label{section:derivation_mqs}
%
%
In this section, the homogenized magnetoquasistatic (MQS) problem is derived. The 
derivation uses two main ingredients: the \emph{MQS assumptions} which makes it 
possible to neglect the displacement currents and the \emph{homogenization} of the 
corresponding multiscale problem. The derivation of this paper is made easier by 
applying the MQS assumptions to the homogenized parabolic hyperbolic (PH) multiscale 
problem that was already carried out in \cite{visintin-tsh-06-b,visintin-tsh-08} 
instead of applying the homogenization theory to the parabolic elliptic (PE) multiscale 
problem derived from the PH multiscale under appropriate assumptions
(see Figure \ref{fig:commutative_diagram_derivation}). 
\added{
In \cite{visintin-tsh-06-b, visintin-tsh-08}, existence and uniqueness of the solution was proved via the 
approximation by time-discretization, the derivation of a priori estimates, 
and the passage to the limit via compensated compactness and compactness by strict convexity.
The homogenized problem was then derived using the two-scale convergence theory
for the fields and the convergence of functionals used to define constitutive laws.
}
In Section \ref{subsection:maxwell-function-spaces} we recall Maxwell's equations
that govern the evolution of electromagnetic fields and we define the function spaces
used for solving these equations in the weak sense. 
In Section \ref{subsection:PH_problem}, we recall the PH multiscale problem and 
its homogenization as done in \cite{visintin-tsh-06-b, visintin-tsh-08}. This 
homogenized problem is then used in Section \ref{subsection:PE_problem} for the 
derivation of the homogenized parabolic elliptic (PE) problem.
In the rest of the section, we use the capital letters P, H end E to denote the 
parabolic, hyperbolic and elliptic problems, respectively. Thus, the PH multiscale 
problem denotes the parabolic hyperbolic multiscale problem whereas the PE--PH 
homogenized problem denotes the homogenized problem with a PE problem at the coarse 
scale and a PH problem are the fine scale. The PE problem corresponds to the MQS 
problem.

\subsection{Maxwell's equations and the function spaces}
\label{subsection:maxwell-function-spaces}

Consider the electromagnetic problem in an open domain $\Omega_T := \Omega \times \mathcal{I}$ 
with $\Omega \subseteq \mathbb{R}^3$ and $\mathcal{I} = (0,T] \subset \mathbb{R}$. 
The electromagnetic fields are governed by the following Maxwell equations and 
constitutive laws \cite{bossavit-modelisation-93, bossavit-computational-98, jackson1999classical}:
\begin{subequations}
    \begin{equation}
        \Curl[]{\bh }  = \bj + \bj_s + \epsilon \partial_t \be, 
        \quad \Curl[]{\be } = - \partial_t \bb, 
        \quad \Div[]{\bb } = 0  \hspace{5mm} \textrm{ in } \Omega \times \mathcal{I},
    \tag{\theequation\,a-c}
    \end{equation}
\label{eq:mqs_maxwell}
\end{subequations}
\vspace{-8mm}
\begin{subequations}
    \begin{equation}
        \bb(\bx, t) = \boldsymbol{\mathcal{B}}(\bh(\bx, t), \bx ), 
        \quad \bj(\bx, t) = \boldsymbol{\mathcal{J}}(\be(\bx, t), \bx ) 
        \hspace{5mm} \forall (\bx, t) \in  \Omega \times \mathcal{I}.
    \tag{\theequation\,a-b}
    \end{equation}
\label{eq:mqs_materiallaw}	
\end{subequations}
\vspace{-4mm}

\noindent The field $\bh$ is the magnetic field, $\bb$ the magnetic flux density, 
$\bj$ the electric current density, $\bj_s$ the imposed electric current density 
(source) and $\be$ the electric field. The material laws \eqref{eq:mqs_materiallaw} 
are expressed in terms of the mappings 
$\boldsymbol{\mathcal{B}}: \mathbb{R}^3 \times \Omega \rightarrow \mathbb{R}^3$ and
$\boldsymbol{\mathcal{J}}: \mathbb{R}^3 \times \Omega \rightarrow \mathbb{R}^3$, 
linear or not, accounting for the magnetic and electric behaviour, respectively. 
The domain $\Omega$ is subdivided into conducting ($\Omega_c$) and nonconducting 
($\Omega_c^C$) parts, the former being where eddy currents can appear.
The boundary of the domain $\Omega$ is denoted $\Gamma$. 
\deleted{
It is the union of 
$\Gamma_e$ and $\Gamma_h$ such that $\Gamma = \Gamma_e \cup \Gamma_h$
and $\Gamma_e \cap \Gamma_h = \emptyset$.
The region $\Gamma_e$ is the part of the boundary where the tangential trace of 
$\be$ (resp. the normal trace of $\bb$) is imposed and $\Gamma_h$ is the part 
of the boundary where the tangential trace of $\bh$ (resp. the normal trace of 
$\bj$) is imposed. 
}
\added{
In Sections \ref{section:mqs} and \ref{section:models} we derive 
the weak solution sof the MQS problem using the magnetic vector potential formulations
\cite{bachinger-mqs-05,jiang-mqs-12,rodriguez-mqs-10,acevedo-mqs-13}. 
In Sections \ref{section:mqs} and \ref{section:models}, some structural restrictions 
on the computational domain are assumed for the existence and the uniqueness of the solution 
\cite{rodriguez-mqs-10,acevedo-mqs-13,bachinger-mqs-05}.
The domain $\Omega$ is assumed to be simply connected with a Lipschitz connected 
boundary $\Gamma$. The conducting domain $\Omega_c$ is an open subset strictly contained 
in $\Omega$ which can be connected or not. In the latter case, $\Omega_c = \cup_{i = 1}^m \Omega_c^i$
where $\Omega_c^i$, i = 1, 2, \ldots, m are connected components of $\Omega_c$. For 
simplicity we assume the non-conducting domain $\Omega_c^C$ to be connected.
The case of a non-connected $\Omega_c^C$ can be also easily treated.
}
The system of equations must further be completed by an initial condition on the 
magnetic flux density assumed to be divergence-free, i.e., $\Div[]{\bb^{0}} = 0$.
The superscript $^0$ is used to denote initial condition, i.e., $\bb^{0} = \bb(\cdot, 0)$. 
This conditions together with (\ref{eq:mqs_maxwell} b) naturally imply Gau\ss \, 
magnetic law (\ref{eq:mqs_maxwell} c). In the rest of this section, we ignore 
Gau\ss \, magnetic law which is automatically fulfilled under Faraday's equation 
(\ref{eq:mqs_maxwell} b) together with this initial condition $\Div[]{\bb^{0}} = 0$ 
(see \cite{visintin-tsh-06-b, visintin-tsh-08}).

%
%
The weak solutions the fullscale, the macroscale and the mesoscale problems must belong 
to the right function spaces. For almost every $t \in \mathcal{I}$, these
functions spaces are defined as the domains of the differential operators $\GradSymb, \CurlSymb$ 
and $\DivSymb$ with appropriate non-homogeneous boundary conditions prescribed on 
the boundary $\Gamma$:
\begin{align}
    \Hone[]{\Omega} &:= \{ u\in\Ltwo{\Omega}: \Grad{u}\in\LLtwo{\Omega}\} , 
    \label{eq:eq_Hone_h}
    \\
    \Hcurl[]{\Omega} &:= \{ \vec{u}\in\LLtwo{\Omega}: \Curl{\vec{u}}\in\LLtwo{\Omega}\} ,
    \label{eq:eq_Hcurl_h}
    \\
    \Hdiv[]{\Omega} &:=\{ \vec{u}\in\LLtwo{\Omega}:\Div{\vec{u}}\in\Ltwo{\Omega}.
\label{eq:eq_Hdiv_h}
\end{align}
The spaces $\Hone[_0]{\Omega}$, $\Hcurl[_0]{\Omega}$, $\Hdiv[_0]{\Omega}$ denote 
the same spaces as the corresponding spaces in \eqref{eq:eq_Hone_h}--\eqref{eq:eq_Hdiv_h} 
with traces equal to zero, i.e.,
\begin{align}
    \Hone[_0]{\Omega} &:= \{ u \in \Hone[]{\Omega}, u|_{\Gamma} = 0\} , 
    \label{eq:eq_Hone_h_zero}
    \\
    \Hcurl[_0]{\Omega} &:= \{ \vec{u} \in \Hcurl[]{\Omega}, \bn \times \bu|_{\Gamma} = \boldsymbol{0}\} ,
    \label{eq:eq_Hcurl_zero}
    \\
    \Hdiv[_0]{\Omega} &:=\{ \vec{u} \in \Hdiv[]{\Omega}, \bn \cdot \bu|_{\Gamma} = 0\}. 
    \label{eq:eq_Hdiv_h_zero}
\end{align} 
The spaces $\Hcurlzero[]{\Omega}$, $\Hdivzero[]{\Omega}$ denote the nullspace
of the operators $\CurlSymb$ and $\DivSymb$, respectively. 
\added{
In Sections \ref{section:mqs} and \ref{section:models} we consider the following 
Bochner spaces for the potentials, solution of the multiscale and the macroscale problems:
}
\begin{equation}
    \LLtwo{0, T; \boldsymbol{V}} \qquad \text{ and } \qquad \LLtwo{0, T; \boldsymbol{V}^{\ast}},
\label{eq:eq_Macro_Space}
\end{equation}
\added{
where $\boldsymbol{V}$ can be any vector space (in Sections \ref{section:mqs} and 
\ref{section:models} we use $\boldsymbol{V} := \Hcurl[_0]{\Omega}$) and $\boldsymbol{V}^{\ast}$ 
is the dual of $\boldsymbol{V}$.
}
The mesoscale problem leads to the solutions that belong to the spaces: 
\vspace{-2mm}
\begin{multline}
    \LLtwo{\mathbb{R}_T^3; \boldsymbol{W}} := \left\{ \bu: \mathbb{R}_T^3 \rightarrow \boldsymbol{W}: \displaystyle{\textcolor{white}{\left(\int_{\mathbb{R}_T^3} \right)^{\frac{1}{2} } } } \right. \\
    \left. \|\bu\|_{\LLtwo{\mathbb{R}_T^3 : \boldsymbol{W}}} := \left(\int_{\mathbb{R}_T^3} \left\|\bu(\bx, t)\right\|_{\boldsymbol{W}}^2 \text{d}t \, \text{d}x \right)^{\frac{1}{2}} < \infty \right\},
\label{eq:eq_Meso_Space}
\vspace{-2mm}
\end{multline}
where the separable Banach space $\boldsymbol{W}$ is defined on the mesoscale domain 
$Y \equiv \Omega_m$. For the homogenized PH problem, two spaces were used in place 
of $\boldsymbol{W}$: the nullspaces $\Hcurlzero[]{\mathcal{Y}}$ and $\Hdivzero[]{\mathcal{Y} }$. 
The symbol $\mathcal{Y}$ is used for functions defined on $Y$ with periodic boundary conditions.

\subsection{Homogenization of the Parabolic Hyperbolic multiscale problem}
\label{subsection:PH_problem}

From now on, we consider $\Omega = \mathbb{R}^3$ and derive the parabolic hyperbolic 
multiscale problem along the lines of \cite{visintin-tsh-06-b,visintin-tsh-08}.
\begin{prob}[\textbf{Parabolic--Hyperbolic (PH) multiscale problem}] 
The PH multiscale problem was derived from Maxwell's equation by neglecting the 
displacement currents with respect to the eddy currents in the conducting domain 
\emph{(}i.e., $\epsilon \partial_t \be^{\varepsilon} \ll \bj^{\varepsilon}$ in $\Omega_c$\emph{)}.
    \begin{subequations}
        \begin{equation}
            \Curl[]{\bh^{\varepsilon} }  = \bj^{\varepsilon} + \bj_s + (1 - \chi_{_{\Omega_c}}) \epsilon^{\varepsilon} \, \partial_t \be^{\varepsilon},
            \qquad \Curl[]{\be^{\varepsilon} } = - \partial_t \bb^{\varepsilon},  
            \tag{\theequation\,a-b}
        \end{equation}
    \label{eq:equations_PH_multiscale}
    \end{subequations}
    \vspace{-8mm}
    \begin{subequations}
        \begin{equation}
            \bb^{\varepsilon}(\bx, t) = \boldsymbol{\mathcal{B}}^{\varepsilon}(\bh^{\varepsilon}(\bx, t), \bx), 
            \, \, \bj^{\varepsilon}(\bx, t) = \boldsymbol{\mathcal{J}}^{\varepsilon} (\be^{\varepsilon}(\bx, t), \bx)
            \, \, \, \, \forall (\bx, t) \in  \mathbb{R}_T^3
            \tag{\theequation\,a-b}
        \end{equation}
    \label{eq:materiallaw_PH_multiscale}	
    \end{subequations}
    \vspace{-0mm}
    \hspace{-3mm} where the function $\chi_{_{\Omega_c}}$ is the characteristic 
    function, different from zero only on the conducting domain $\Omega_c$. The 
    superscript $^{\varepsilon}$ is used to denote the multiscale dependency of 
    the fields. All derivatives are defined in the distribution sense.
    \label{prob:PH_multiscale_problem}
\end{prob}
In \cite{visintin-tsh-06-b,visintin-tsh-08}, Gau\ss \, magnetic law $\Div[]{\bb^{\varepsilon}} = 0$ 
was ensured by imposing the initial condition on $\bb^{\varepsilon 0}$ such that $\Div[]{\bb^{\varepsilon 0} } = 0$. 
The material laws \eqref{eq:mqs_materiallaw} are expressed in terms of the mappings 
$\boldsymbol{\mathcal{B}}^{\varepsilon} : \mathbb{R}^3 \times \Omega \rightarrow \mathbb{R}^3$ 
and $\boldsymbol{\mathcal{J}}^{\varepsilon}  : \mathbb{R}^3 \times \Omega \rightarrow \mathbb{R}^3$
defined by:
\begin{equation}
    \boldsymbol{\mathcal{B}}^{\varepsilon}(\bh^{\varepsilon}, \bx) = \boldsymbol{\mathcal{\bar{B}}}(\bh^{\varepsilon}, \bx, \bx/\varepsilon),
    \qquad \boldsymbol{\mathcal{J}}^{\varepsilon}(\be^{\varepsilon}, \bx) = \boldsymbol{\mathcal{\bar{J}}}(\be^{\varepsilon}, \bx, \bx/\varepsilon),
\label{eq:materiallaw_multiscale_definition}
\end{equation}
where the operators $\boldsymbol{\mathcal{\bar{B}}} : \mathbb{R}^3 \times \Omega \times Y \rightarrow \mathbb{R}^3$
and  $\boldsymbol{\mathcal{\bar{J}}} : \mathbb{R}^3 \times \Omega \times \times Y \rightarrow \mathbb{R}^3$
are used to represent two-scale composite materials for which the characteristic 
length at the mesoscale is $\varepsilon$. By abuse of notation, we use $\boldsymbol{\mathcal{B}}$ 
and $\boldsymbol{\mathcal{J}}$ instead of $\boldsymbol{\mathcal{\bar{B}}}$
and $\boldsymbol{\mathcal{\bar{J}}}$ in the rest of the text.
For the analytical and theoretical study of the multiscale Problem \ref{prob:PH_multiscale_problem} 
we assume that the nonlinear mapping $\boldsymbol{\mathcal{B}}$ is maximal monotone 
and therefore it can be derived by the minimization of a convex, lower-semicontinous 
functional. It also has an inverse ${\boldsymbol{\mathcal{B}}^{-1} \equiv \boldsymbol{\mathcal{H}}}$ 
that can be derived from a conjuguate convex, lower semi-continuous functional
\cite{ekeland-convexanalysis-69, evans-pde-10, rockafellar-convexanalysis-69}.
This covers cases of linear and nonlinear reversible magnetic laws. However, one
of the major advantages of the computational homogenization approach proposed in
Section \ref{section:models} is the inclusion of hysteretic laws in the
numerical model by means of classical hysteresis models (e.g. Preisach,
Jiles-Atherton, etc.). We will thus lift this hypothesis once we consider the
computational framework.  We will still assume that the mapping $\boldsymbol{\mathcal{J}}$ 
is maximal monotone and has an inverse $\boldsymbol{\mathcal{J}}^{-1} \equiv 
\boldsymbol{\mathcal{E}}$. In practice, this assumption holds as the materials we 
consider in this paper are electrically linear.

Problem \ref{prob:PH_multiscale_problem} has been extensively analyzed. A homogenized 
problem with \emph{coarse} and \emph{fine} problems was derived considering 
some assumptions on the constitutive laws, the initial conditions (IC) and the 
current source $\bj_s$. These assumptions are recalled in \textsc{Assumptions} 
\ref{ass:PH_multiscale_IC_Source}--\ref{ass:PH_multiscale_Convergence_IC}
\begin{ass}[\textbf{Regularity of the IC and the sources}]
Assume that the \\initial conditions $\bb^{\varepsilon 0}$ and $\be^{\varepsilon 0}$ 
and the source $\bj_s$ fulfill the following regularity conditions:
\vspace{-2mm}
    \begin{subequations}
        \begin{equation}
            \hspace{-2.5mm} \bb^{\varepsilon 0} \in \LLtwo{\mathbb{R}^3},
            \, \, \, \be^{\varepsilon 0} \in \LLtwo{\Omega_c^C},
            \, \, \, \bj_s \in \LLtwo{\Omega_s \times \mathcal{I}},
            \, \, \, \Div[]{\bb^{\varepsilon 0} } = 0,
            \, \, \, \Div[]{\bj_s } = 0.
            \tag{\theequation\,a-e}
        \end{equation}
        \label{eq:equations_PH_multiscale_IC_Source}
    \end{subequations}
    \label{ass:PH_multiscale_IC_Source}
\end{ass}
\vspace{-6mm}
Equation (\ref{eq:equations_PH_multiscale_IC_Source} d) together with (\ref{eq:equations_PH_multiscale} b) 
ensures Gau\ss \, magnetic law $\Div[]{\bb^{\varepsilon}} = 0$.  
\vspace{-0mm}
\begin{ass}[\textbf{Assumptions on the constitutive laws}]%
    Assume that the \\electrical law is given by $\bj^{\varepsilon} = \sigma^{\varepsilon} \, 
    \be^{\varepsilon}$ where the electrical conductivity $\sigma^{\varepsilon}$ is 
    definite positive in $\Omega_c$ and that the mapping $\boldsymbol{\mathcal{B}}$ 
    is maximal monotone.
    \label{ass:PH_multiscale_ConstitutiveLaws}
\end{ass}
\vspace{-0mm}
These restrictions on the mappings cover a wide range of material laws usually 
encountered in applications. They cover the linear electrical materials, the linear and 
the nonlinear reversible magnetic materials as well as soft magnetic materials for which 
the hysteresis loop can be approximated using the maximal monotone operators. However,
the hard magnetic materials are not covered.
\begin{ass}[\textbf{Convergence of the initial conditions}]%
    Assume that the \\initial conditions $\bb_0^{\varepsilon}$ and $\be_0^{\varepsilon}$ 
    converge in the classical and the two-scale senses, i.e.:
    \vspace{-2mm}
    \begin{subequations}
        \begin{equation}
            \begin{aligned}
                &\bb^{\varepsilon 0} \underset{2}{\rightharpoonup} \bb_0^{0}\, \, \emph{\text{ in }} \, \, \LLtwo{\mathbb{R}^3 \times \mathcal{Y}}, 
                \qquad &&\bb^{\varepsilon 0} \rightharpoonup \left< \bb_0^{0} \right>_{Y} = \bb_M^{0} \, \, \emph{\text{ in }} \, \, \LLtwo{\mathbb{R}^3},
            \end{aligned}
        \tag{\theequation\,a-b}
        \end{equation}  
    \end{subequations}
    \vspace{-8mm}      
    \begin{subequations}
        \begin{equation}
            \begin{aligned}
                &\be^{\varepsilon 0} \underset{2}{\rightharpoonup} \be_0^{0} \, \, \emph{\text{ in }} \, \, \LLtwo{\Omega_c^C \times \mathcal{Y}},
                \qquad &&\be^{\varepsilon 0} \rightharpoonup \left< \be_0^{0} \right>_{Y} = \be_M^{0} \, \, \emph{\text{ in }} \, \, \LLtwo{\Omega_c^C}.
            \end{aligned}
        \tag{\theequation\,a-b}
        \end{equation}
        \label{eq:equations_PH_multiscale_Convergence_IC}
    \end{subequations}
    \vspace{-0mm}
    \hspace{-1.50mm}These fields are used as initial conditions for the fine and 
    the coarse problem, respectively. The curly brackets $\left< \bff \right>_{Y}$ 
    are used to denote the average of the function $\bff$ over the cell domain $Y$,
    i.e., $$\left< \bff \right>_{Y} = \frac{1}{|Y|} \int_{Y} \bff \text{d} y = \frac{1}{|\Omega_{m}|} \int_{\Omega_{m}} \bff \text{d} y = \left< \bff \right>_{\Omega_{m}}$$ where
    $|Y|$ is used to denote the volume of the domain $Y \equiv \Omega_m$.
\label{ass:PH_multiscale_Convergence_IC}
\end{ass}

\vspace{5mm}

Using \textsc{Assumption} \ref{ass:PH_multiscale_IC_Source} for the IC and the 
source term and \textsc{Assumption} \ref{ass:PH_multiscale_ConstitutiveLaws} for the constitutive 
laws, the following PH--PH homogenized problem was derived from the multiscale Problem 
\ref{prob:PH_multiscale_problem} (\cite{visintin-tsh-08}): 
\begin{prob}[\textbf{PH--PH homogenized problem}] 
The PH--PH homogenized \\problem has been derived from Problem \ref{prob:PH_multiscale_problem} 
with the following two coarse and fine problems\emph{:}\\
\hspace{-4mm}
\emph{\textbf{Coarse problem:} find $\bh_M, \be_M, \bb_M, \bj_M \in \LLtwo{\mathbb{R}_T^3}$ such that} 
    \begin{subequations}
        \begin{equation}
            \Curl[_x]{\bh_M} = \bj_M + \bj_s + (1 - \chi_{_{\Omega_c}}) \epsilon_M \partial_t \be_M, 
            \, \, \, \,  \Curl[_x]{\be_M} = - \partial_t \bb_M,
            \tag{\theequation\,a-b}
        \end{equation}
    \label{eq:equations_PH_PH_homogenized_macro}
    \end{subequations}
    \vspace{-8mm}
    \begin{subequations}
        \begin{equation}
            \bb_M = \boldsymbol{\mathcal{B}}_M(\bh_M, \bx), 
            \, \, \bj_M = \boldsymbol{\mathcal{J}}_M(\be_M, \bx),
            \, \, \text{ for a.e. } (\bx, t) \in  \mathbb{R}^3 \times \mathcal{I}. 
            \tag{\theequation\,a-b}
        \end{equation}
    \label{eq:materiallaw_PH_PH_homogenized_macro}	
    \end{subequations}
    \vspace{-4mm}
    \\
\emph{\textbf{Fine problem:} find $\bh_0, \be_0 \in \LLtwo{\mathbb{R}_T^3\!\!:\!\!\Hcurlzero[]{\mathcal{Y} } }$
and $\bh_1, \be_1, \bb_0, \bj_0 \in \LLtwo{\mathbb{R}_T^3\!\!:\!\!\Hdivzero[]{\mathcal{Y} } } $ such that}
    \begin{subequations}
        \begin{multline}
            \Curl[_x]{\bh_M} + \Curl[_y]{\bh_1}  = \bj_0 + (1 - \chi_{_{\Omega_c}}) \epsilon \partial_t \be_0,
            \\ 
            \Curl[_x]{\be_M} + \Curl[_y]{\be_1} = - \partial_t \bb_0, \textcolor{white}{Filling text }
            \tag{\theequation\,a-b}
        \end{multline}
    \label{eq:equations_PH_PH_homogenized_meso}
    \end{subequations}
    \vspace{-8mm}
    \begin{subequations}
        \begin{equation}
            \bb_0 = \boldsymbol{\mathcal{B}}(\bh_0, \bx, \by), 
            \, \, \bj_0 = \boldsymbol{\mathcal{J}}(\be_0, \bx, \by),
            \, \, \text{ for a.e. } (\bx, \by, t) \in  \mathbb{R}^3 \times Y \times \mathcal{I}_t. 
            \tag{\theequation\,a-b}
        \end{equation}
    \label{eq:materiallaw_PH_PH_homogenized_meso}	
    \end{subequations}
    The macroscale fields are obtained as averages of the zero order terms, i.e., 
    $\bff_M = \langle {\bff}_0 \rangle_Y$. All the derivatives are defined in the 
    distribution sense. 
    \label{prob:PH_PH_homogenized_problem}
\end{prob}
Equation (\ref{eq:equations_PH_PH_homogenized_macro} b) together with $\Div[_x]{\bb_M^0} = 0$ 
imply the coarse scale Gau\ss \, magnetic law $\Div[_x]{\bb_M} = 0$. The equations 
of the fine scale (\ref{eq:equations_PH_PH_homogenized_meso} a-b)--(\ref{eq:materiallaw_PH_PH_homogenized_meso} a-b)
involve the nullspaces that can be decomposed 
as \cite{visintin-tsh-08,visintin-tsh-11,pankov-gconv-97,monk2003finite}:
\begin{align}
    &\Hcurlzero[]{\mathcal{Y} } = \mathbb{R}^3 \, \oplus \, \Hcurlzero[_{*}]{\mathcal{Y} } = \mathbb{R}^3 \, \oplus \, \GradSymb_y \, \Hone[_*]{\mathcal{Y}},
    \label{eq:space-decomposition_Hcurl}\\
    &\Hdivzero[]{\mathcal{Y} } = \mathbb{R}^3 \, \oplus \, \Hdivzero[_{*}]{\mathcal{Y} } = \mathbb{R}^3 \, \oplus \, \CurlSymb_y \, \Hcurl[_*]{\mathcal{Y}}.
    \label{eq:space-decomposition_Hdiv}
\end{align}
Using the decompositions in \eqref{eq:space-decomposition_Hcurl} and \eqref{eq:space-decomposition_Hdiv}, 
each field $\bff_0$ of $\Hcurlzero[]{\mathcal{Y} }$ or $\Hdivzero[]{\mathcal{Y} }$ 
can be written as the sum of an average value $\left< \bff_0 \right>_Y \in \mathbb{R}^3$ 
and a zero average perturbation $\tilde{\bff}_0$. 
The second equalities in \eqref{eq:space-decomposition_Hcurl} and \eqref{eq:space-decomposition_Hdiv}
are obtained using the Helmholtz decomposition of $\LLtwo[_*]{ \mathcal{Y} }$: 
\begin{equation}
    \LLtwo[_*]{ \mathcal{Y} } = \GradSymb_y \, \Hone[_*]{\mathcal{Y}} \, \oplus \, \CurlSymb_y \, \Hcurl[_*]{\mathcal{Y}}
\label{eq:space-decomposition_L2}
\end{equation} 
which applies for fields with periodic boundary conditions. Indeed, the subspace 
of gradients of a harmonic function which appears in the general decomposition of 
$\boldsymbol{L}^2$ fields is dismissed in the case of periodic functions 
\eqref{eq:space-decomposition_L2} and for $\Omega = \mathbb{R}^n$ 
(\cite{deriaz-helmholtzdecomposition-09, harouna-helmholtzdecomposition-10}).

The decomposition \eqref{eq:space-decomposition_L2} was used by Visintin for the convergence 
of functionals used to derive the nonlinear magnetic material laws.
For almost every $(\bx, t) \in \mathbb{R}_T^3$, the decompositions in 
\eqref{eq:space-decomposition_Hcurl}--\eqref{eq:space-decomposition_Hdiv} 
leads to the decompositions of the first order terms $\be_0 = \be_M + \Grad[_y]{v_c}$ and 
$\bb_0 = \bb_M + \Curl[_y]{\ba_c}$ with $v_c \in \Hone[_*]{\mathcal{Y}}$ and 
$\ba_c \in \Hcurl[_*]{\mathcal{Y}}$.

If the mappings $\boldsymbol{\mathcal{B}}$ and $\boldsymbol{\mathcal{J}}$ are maximal 
monotone then the mappings $\boldsymbol{\mathcal{B}}_M$ and $\boldsymbol{\mathcal{J}}_M$ 
are also maximal monotone and derived using the approach described in the following paragraphs. 
Their inverses $\boldsymbol{\mathcal{H}}_M \equiv \boldsymbol{\mathcal{B}}_M^{-1}$
and $\boldsymbol{\mathcal{E}}_M \equiv \boldsymbol{\mathcal{J}}_M^{-1}$ can therefore 
be determined by minimizing the convex conjugate functionals and determined by means 
of the mesoscale problems hereafter \cite{visintin-tsh-08}. 

\noindent For the mapping $\boldsymbol{\mathcal{H}}_M$: find $\ba_c \in \Hcurl[_{*}]{ \mathcal{Y} }$ such that
\begin{equation}
    \left(\boldsymbol{\mathcal{H}}\left(\bb_M + \Curl[_y]{\ba_c }, \bx, \by \right), \Curl[_y]{\ba_c'} \right) = 0\,, \qquad \forall \ba_c' \in \Hcurl[_{*}]{ \mathcal{Y} }
\label{eq:b_conform_macro_materiallaw_h}
\end{equation}
and then derive: $\boldsymbol{\mathcal{H}}_M (\bb_M + \bB_c, \bx) = \langle \boldsymbol{\mathcal{H}} (\bb_M + \bB_c, \bx, \by ) \rangle_Y.$

\noindent For the mapping $\boldsymbol{\mathcal{J}}_M$: find $v_c \in \Hone[_{*}]{\mathcal{Y} }$ such that
\begin{equation}
    \left(\boldsymbol{\mathcal{J}}\left(\be_M + \Grad[_y]{v_c }, \bx, \by \right), \Grad[_y]{v_c'}\right) = 0\, \qquad \forall v_c' \in \Hone[_{*}]{ \mathcal{Y} },
    \label{eq:b_conform_macro_materiallaw_j}
\end{equation}
and then derive: $\boldsymbol{\mathcal{J}}_M (\be_M + \bE_c, \bx) 
= \langle \boldsymbol{\mathcal{J}} (\be_M + \bE_c, \bx\, by) \rangle_Y.$

The operators 
\begin{align}
    &\bB_c: Y \times \mathbb{R}^3 \rightarrow \LLtwo[_*]{ \mathcal{Y} }: (\by, \bb_M)  \mapsto \bb_c = \bB_c(\by, \bb_M)
    \label{eq:solution-operators-b}
    \\
    &\bE_c: Y \times \mathbb{R}^3 \rightarrow \LLtwo[_*]{\mathcal{Y} }: (\by, \be_M)  \mapsto \be_c = \bE_c(\by, \be_M)
    \label{eq:solution-operators-e}
\end{align} 
are solution operators for the mesoscale problems with 
$\bb_c = \bB_c(\by, \bb_M)=\Curl[_y]{\ba_c}$ and $\be_c = \bE_c(\by, \be_M)=\Grad[_y]{\bv_c}$.
If the mappings $\boldsymbol{\mathcal{H}}$ and $\boldsymbol{\mathcal{J}}$ are 
linear, problem \eqref{eq:b_conform_macro_materiallaw_h}--\eqref{eq:b_conform_macro_materiallaw_j} 
is equivalent to the cell problem obtained using the asymptotic expansion theory 
\cite{sanchez-homogenization-87, bensoussan-ahm-11, visintin-tsh-11}. The dual 
formulation allows to define similar problems for the constituitive laws 
$\boldsymbol{\mathcal{B}}_M \equiv \boldsymbol{\mathcal{H}}_M^{-1}$ and 
$\boldsymbol{\mathcal{E}}_M = \boldsymbol{\mathcal{J}}_M^{-1}$.
%
%
\subsection{Homogenization of the Parabolic Elliptic multiscale problem}
\label{subsection:PE_problem}

The MQS problem can be derived by applying the MQS assumption to Maxwell's equations. 
This assumption can be derived by comparing the following physical parameters of 
the problem: $L_c$ and $L_f$ which are the coarse and fine scale characteristic 
lengths (e.g., the sizes of the coarse and the fine domains), $\lambda_f$ and $\lambda_M$ 
which are the coarse and the fine wavelengths respectively and $\delta_c $ and $\delta_f$, 
the coarse and the fine skin depths, respectively. The wavelengths 
$\lambda_f = 2 \pi/(\omega \, \sqrt{\mu \epsilon})$ and $\lambda_M = 2 \pi/(\omega \, \sqrt{\mu_M \epsilon_M})$ 
and the skin depths are defined by $\delta_f = \sqrt{2/\omega \, \sigma \, \mu}$
and $\delta_M = \sqrt{2/\omega \, \sigma_M \, \mu_M}$ where $\sigma_M$ and $\epsilon_M$ are 
the homogenized electric conductivity and permittivity that can be obtained by 
solving a linear electrokinetic and electrostatic cell problems \cite{niyonzima-chm-13, niyonzima-chm-13-1} 
and $\mu_M$ is the nonlinear homogenized magnetic permeability which can be determined 
from \eqref{eq:b_conform_macro_materiallaw_h}. Additionally, the magnetostatic problem
can de derived from the MQS problem by neglecting the eddy currents if the MS 
\textsc{Assumption} \ref{ass:PE_EddyCurrents_Assumptions} is fulfilled.
The conditions that lead to the magnetoquasistatic and the magnetostatic problems are stated in
\textsc{Assumptions} \ref{ass:PE_MQS_Assumptions}--\ref{ass:PE_EddyCurrents_Assumptions}.
\begin{ass}[\textbf{MQS assumption}]
Displacement currents at the coarse and fine scales can be neglected if the following 
conditions are fulfilled:
    \begin{enumerate}
        \item The displacement currents at the coarse scale $(1 - \chi_{_{\Omega_c}}) \epsilon_M \partial_t \be_M$ 
            can be neglected if $\lambda_c/L_c \gg 1$.    
        \item The displacement currents at the fine scale $(1 - \chi_{_{\Omega_c}}) \epsilon \partial_t \be_0$ 
            can be neglected if $\lambda_f/L_f \gg 1$.
    \end{enumerate}
    \label{ass:PE_MQS_Assumptions}
\end{ass}
\begin{ass}[\textbf{MS assumption}]
The coarse-scale and the fine-scale eddy currents can be neglected if the following 
conditions are fulfilled:
    \begin{enumerate}
        \item The coarse scale eddy currents $\bj_M$ can be neglected if there is no 
        net coarse scale eddy currents \emph{(}e.g.: in the case of perfect insulation\emph{)} 
        or if $\delta_c/L_c \gg 1$.
        \item The mesoscale eddy currents $\bj_0$ can be neglected if there are no 
        conducting materials in the cell unit (i.e., $\Omega_{mc} = \emptyset$) or 
        if $\delta_f/L_f \gg 1$.
    \end{enumerate}
    \label{ass:PE_EddyCurrents_Assumptions}
\end{ass}
\begin{table}[h!]\caption{Type of problems depending on the predefined physical parameters of the problem.}
\label{tab:list-derived-problems}
\medskip\centering
\begin{tabular}{| c | c | c | c | c | c | c | c |}
\hline
\# Problem & $\lambda_c/L_c$ & $\lambda_f/L_f$ & $\delta_c/L_c$ & $\delta_f/L_f$ & Multiscale & Coarse     & Fine \\
\hline
 (1)   &$\simeq 1$       & $\simeq 1$      & $\simeq 1$     & $\simeq 1$     & PH         & PH         & PH   \\
 (2)   &$\gg 1$          & $\simeq 1$      & $\simeq 1$     & $\simeq 1$     & PH         & P          & PH   \\
 (3)   &$\simeq 1$       & $\simeq 1$      & $\gg 1$        & $\simeq 1$     & PH         & H          & PH   \\
 (4)   &$\simeq 1$       & $\simeq 1$      & $\gg 1$        & $\gg 1$        & H          & H          & H    \\
 (5)   &$\gg 1$          & $\gg 1$         & $\simeq 1$     & $\simeq 1$     & PE         & PE         & PE   \\
 (6)   &$\gg 1$          & $\gg 1$         & $\gg 1$        & $\simeq 1$     & PE         & E          & PE   \\
 (7)   &$\gg 1$          & $\gg 1$         & $\simeq 1$     & $\gg 1$        & PE         & PE         & E    \\
 (8)   &$\gg 1$          & $\gg 1$         & $\gg 1$        & $\gg 1$        & E          & E          & E    \\
\hline
\end{tabular}
\end{table}
The combination of the parameters defined above lead to the multiscale and homogenized
problems defined in \textsc{Table} \ref{tab:list-derived-problems}. 
In this paper we focus on the PE multiscale problem \ref{prob:PE_multiscale_problem} 
derived using \textsc{Assumption} \ref{ass:PE_MQS_Assumptions}.
\begin{prob}[\textbf{Parabolic--Elliptic (PE) multiscale problem}] 
This problem can be derived from Problem \ref{prob:PH_multiscale_problem} if point 
2 of Assumption \ref{ass:PE_MQS_Assumptions} is fulfilled. In that case, the displacement currents 
$\epsilon \partial_t \be^{\varepsilon}$ can be neglected in the entire domain leading 
to the following equations: 
\vspace{-2mm}
\begin{subequations}
        \begin{equation}
            \Curl[]{\bh^{\varepsilon} }  = \bj^{\varepsilon} + \bj_s,
            \qquad \Curl[]{\be^{\varepsilon} } = - \partial_t \bb^{\varepsilon},  
            \tag{\theequation\,a-b}
        \end{equation}
\label{eq:equations_PE_multiscale}
\end{subequations}
\vspace{-6mm}
\begin{subequations}
    \begin{equation}
        \bb^{\varepsilon} = \boldsymbol{\mathcal{B}}^{\varepsilon}(\bh^{\varepsilon}, \bx), 
        \qquad \bj^{\varepsilon} = \boldsymbol{\mathcal{J}}^{\varepsilon} (\be^{\varepsilon}, \bx)
        \qquad \text{ for a.e. } (\bx, t) \in  \mathbb{R}_T^3.
        \tag{\theequation\,a-b}
    \end{equation}
\label{eq:materiallaw_PE_multiscale}	
\end{subequations}
Gau\ss \, magnetic law $\Div[]{\bb^{\varepsilon} } = 0$ is automatically verified 
if the initial condition $\Div[]{\bb_0^{\varepsilon} } = 0$ is imposed.
\label{prob:PE_multiscale_problem}
\end{prob}
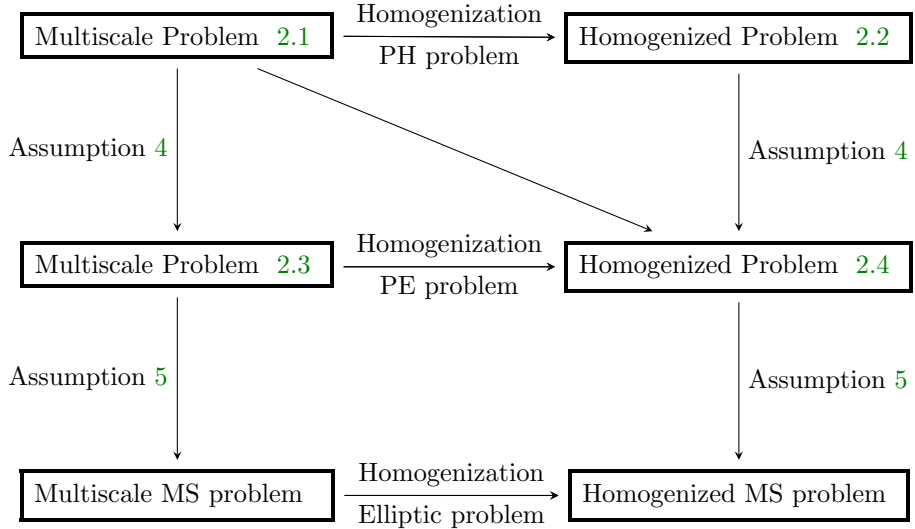
\begin{figure}[h!]
    \centering
    \begin{tikzpicture}[scale=2.750]
        \matrix (m) [matrix of math nodes,row sep=6em,column sep=8em,minimum width=4em]
        {
            \mybox{\text{{Multiscale Problem }} \ref{prob:PH_multiscale_problem} } & \mybox{\text{{Homogenized Problem }} \ref{prob:PH_PH_homogenized_problem} } \\
            \mybox{\text{{Multiscale Problem }} \ref{prob:PE_multiscale_problem} } & \mybox{\text{{Homogenized Problem }} \ref{prob:P_P_homogenized_problem} } \\
            \mybox{\text{{Multiscale MS problem }} } & \mybox{\text{{Homogenized MS problem}} } \\
        };
        \path[-stealth]
        (m-1-1) edge node [left] {Assumption \ref{ass:PE_MQS_Assumptions}} (m-2-1)
                edge node [above] {Homogenization} (m-1-2)
                edge node [below] {PH problem} (m-1-2)
                edge node [above] {} (m-2-2)
        (m-2-1.east|-m-2-2) edge node [above] {Homogenization} (m-2-2)
        (m-2-1.east|-m-2-2) edge node [below] {PE problem} (m-2-2)
        (m-1-2) edge node [right] {Assumption \ref{ass:PE_MQS_Assumptions}} (m-2-2)
        (m-2-1) edge node [left] {Assumption \ref{ass:PE_EddyCurrents_Assumptions}} (m-3-1)

        (m-3-1.east|-m-3-2) edge node [above] {Homogenization} (m-3-2)
        (m-3-1.east|-m-3-2) edge node [below] {Elliptic problem} (m-3-2)
        (m-2-2) edge node [right] {Assumption \ref{ass:PE_EddyCurrents_Assumptions}} (m-3-2);
    \end{tikzpicture}
    \caption{Diagramm illustrating the derivation of the homogenized MQS and magnetostatic 
    problem. The notation ``\textnormal{MS}'' in the diagram stands for magnetostatic.}
    \label{fig:commutative_diagram_derivation}
\end{figure}

The homogenized PE problem can be derived from the PE multiscale problem \ref{prob:PE_multiscale_problem}
using the two-scale and the convergence of functionals as done in \cite{visintin-tsh-06-b,visintin-tsh-08}.
This approach was used in \cite{niyonzima-chm-16} where the multiscale Problem \ref{prob:PE_multiscale_problem}
was solved using the vector potential formulation and then homogenized. In this paper we choose a different 
approach. We use results of the homogenized PH problem and apply the MQS \textsc{Assumption}
\ref{ass:PE_MQS_Assumptions} to derive the homogenized PE problem as illustrated 
in the commutative diagram in Fig. \ref{fig:commutative_diagram_derivation}.
If points 1. and 2. of Assumption \ref{ass:PE_MQS_Assumptions} are valid, the coarse-scale
and the fine-scale displacement currents can be neglected leading to the following 
PE homogenized problem.

\begin{prob}[\textbf{PE--PE homogenized problem}] 
This problem can be derived from the multiscale Problem \ref{prob:PH_PH_homogenized_problem} with
the following coarse and fine problems:\\
\vspace{-0mm}
\hspace{-2mm}\emph{\textbf{Coarse problem:} find $\bh_M, \be_M, \bb_M, \bj_M \in \LLtwo{\mathbb{R}_T^3}$ such that}
    \begin{subequations}
        \begin{equation}
            \Curl[_x]{\bh_M} = \bj_M + \bj_s, 
            \, \, \, \,  \Curl[_x]{\be_M} = - \partial_t \bb_M,
            \tag{\theequation\,a-b}
        \end{equation}
    \label{eq:equations_P_P_homogenized_macro}
    \end{subequations}
    \vspace{-7mm}
    \begin{subequations}
        \begin{equation}
            \bb_M = \boldsymbol{\mathcal{B}}_M(\bh_M, \bx), 
            \, \, \bj_M = \boldsymbol{\mathcal{J}}_M(\be_M, \bx),
            \, \, \text{ for a.e. } (\bx, t) \in  \mathbb{R}^3 \times \mathcal{I}. 
            \tag{\theequation\,a-b}
        \end{equation}
    \label{eq:materiallaw_P_P_homogenized_macro}	
    \end{subequations}
    \vspace{-5mm}
    \\
\emph{\textbf{Fine problem:} find $\bh_0, \be_0 \in \LLtwo{\mathbb{R}_T^3\!\!:\!\!\Hcurlzero[]{\mathcal{Y} } }$
and $\bh_1, \be_1, \bb_0, \bj_0 \in \LLtwo{\mathbb{R}_T^3\!\!:\!\!\Hdivzero[]{\mathcal{Y} } } $ such that}
    \begin{subequations}
        \begin{equation}
            \Curl[_x]{\bh_M} + \Curl[_y]{\bh_1}  = \bj_0,
            \qquad \Curl[_x]{\be_M} + \Curl[_y]{\be_1} = - \partial_t \bb_0,
            \tag{\theequation\,a-b}
        \end{equation}
    \label{eq:equations_P_P_homogenized_meso}
    \end{subequations}
    \vspace{-8mm}
    \begin{subequations}
        \begin{equation}
            \bb_0 = \boldsymbol{\mathcal{B}}(\bh_0, \bx, \by), 
            \, \, \bj_0 = \boldsymbol{\mathcal{J}}(\be_0, \bx, \by),
            \, \, \text{ for a.e. } (\bx, \by, t) \in  \mathbb{R}^3 \times Y \times \mathcal{I}. 
            \tag{\theequation\,a-b}
        \end{equation}
    \label{eq:materiallaw_P_P_homogenized_meso}	
    \end{subequations}
    Equations (\ref{eq:equations_P_P_homogenized_macro} a-b) and 
    (\ref{eq:equations_P_P_homogenized_meso} a-b) are defined in the distribution sense.
    \label{prob:P_P_homogenized_problem}
\end{prob}
\vspace{10mm}
%
%
\section{The magnetoquasistatic approximation}
\label{section:mqs}

In this section we develop the weak formulations for the multiscale problem 
(\ref{eq:equations_PE_multiscale} a-b)--(\ref{eq:materiallaw_PE_multiscale} a-b). 
We omit the superscript $^{\varepsilon}$ to lighten the contents of the section.

\subsection{Magnetic flux density conforming formulations: dynamic case}
\label{subsection:mqs_b_dynamic}
We assume the electrical constitutive law in (\ref{eq:mqs_materiallaw} b) to be 
of the form $\bj = \sigma \be$ where $\sigma$ is the electric conductivity assumed to be piecewise constant.  
We want to solve (\ref{eq:equations_PE_multiscale} a-b)--(\ref{eq:materiallaw_PE_multiscale} a-b) using
the so-called magnetic flux density conforming formulation
\cite{bossavit-computational-98,00.compumag99.dular.MagDynHA,ren-hybrid-90}. 

From Gau\ss \, magnetic law $\Div[]{\bb} = 0$ and (\ref{eq:equations_PE_multiscale} b), 
the electric field $\be$ and the magnetic flux density $\bb$ can be expressed in 
terms of the so-called modified magnetic vector potential $\ba$ as 
\begin{equation}
    \bb = \Curl[]{\ba} \qquad \text{and} \qquad \be = -\partial_t \ba.
    \label{eq:mqs_b_conform_potentials}
\end{equation}
We therefore derive the following weak form of Amp\`{e}re's equation
(\ref{eq:equations_PE_multiscale} a) (see \cite{bachinger-mqs-05,jiang-mqs-12}):\\ 
find $\ba \in \Ltwo{0, T; \boldsymbol{V}}$ with $\partial_t \ba \in \Ltwo{0, T; \boldsymbol{V}^{\ast}}$
such that
\begin{equation}
    \ivol[_\Omega]{\bh }{\Curl{{\ba}'}} 
    - \ivol[_{\Omega}]{{\bj}}{{\ba}'}
    =
    \ivol[_{\Omega_s}]{\vec{j}_s}{\vec{a}'},
    \label{eq:mqs_b_conform_ampere_weakform}
\end{equation}
holds for ${\ba}' \in \boldsymbol{V}$. The vector potential $\ba$ is not uniquely 
defined and a gauge condition must be imposed \cite{bachinger-mqs-05,ledger-formulations-10}.
The space $\boldsymbol{V} = \Hcurl[_0]{\Omega}$ with the homogeneous boundary conditions
has been defined in \eqref{eq:eq_Macro_Space} and its use leads to the neglect of the boundary 
term $\isur[_{\Gamma}]{\bn \times \vec{h}}{\vec{a}'}$ in \eqref{eq:mqs_b_conform_ampere_weakform}.

The magnetic vector potential formulation for the three-dimensional  MQS problem leads to the following problem.
\begin{prob}[\textbf{Weak form of the three-dimensional  MQS problem}] Using (\ref{eq:materiallaw_PE_multiscale} a-b) 
    and introducing \eqref{eq:mqs_b_conform_potentials} in \eqref{eq:mqs_b_conform_ampere_weakform}, 
    one gets the weak form: find $\ba \in \Ltwo{0, T; \boldsymbol{V}}$ with $\partial_t \ba \in \Ltwo{0, T; \boldsymbol{V}^{\ast}}$ such that
    \begin{equation}
        \ivol[_{\Omega_c}]{\sigma\,\partial_t \vec{a}}{\vec{a}'}    
        + \ivol[_\Omega]{\bh }{\Curl{\vec{a}'}} 
        =
        \ivol[_{\Omega_s}]{\vec{j}_s}{\vec{a}'},
        \label{eq:mqs_b_conform_magdyn_a}
    \end{equation}
    for all $\vec{a}' \in \boldsymbol{V}$.
\end{prob}
\deleted{
The boundary term in~\eqref{eq:mqs_b_conform_magdyn_a} contains the 
tangential component of the magnetic field which is subject to the natural boundary 
condition on $\Gamma_h$. In this paper we only consider the case of homogeneous 
Neumann condition. In practice, such Neumann condition can be used for imposing 
the symmetry on a plane crossed by a zero electric current or for imposing homogeneous
boundary conditions on a perfect magnetic material, i.e., a material with the magnetic
permeability $\mu \sim \infty$.
}

The two-dimensional case with all currents perpendicular to the section is obtained 
by assuming the source current density $\bj_s = j_s(x, y) \boldsymbol{1}_z$ where 
$\boldsymbol{1}_z$ is the unit vector along the $z$ axis. If the electric conductivity 
$\sigma$ is such that $\sigma_{13} = 0 =\sigma_{23}$, then $z$-components of
the magnetic field $\bh$ and of the magnetic flux density $\bb$ vanish and
it is possible to derive the magnetic flux density $\vec{b}$ from a scalar
potential $a_z(\bx, \by)$ with $\vec{a} = a_z \boldsymbol{1}_z$. In this
case the $\CurlSymb$ operator can be expressed in terms of the $\GradSymb$
operator as $\CurlSymb := \boldsymbol{1}_z \times \GradSymb$ and the
magnetic flux density reads $\vec{b} = \Curl[]{\vec{a}} = \boldsymbol{1}_z
\times \Grad[]{a_z}$. The weak form of the two-dimensional  problem can be derived from \eqref{eq:mqs_b_conform_magdyn_a}.
\begin{prob}[\textbf{Weak form of the two-dimensional  MQS problem}] The weak form of the magnetic 
vector potential formulation of a two-dimensional  MQS problem reads:
    find $a_z \in \Ltwo{0, T; \Hone[_0]{\Omega} }$ with $\partial_t a_z \in \Ltwo{0, T; H^{-1}(\Omega) }$ such that
    \begin{equation}
        \ivol[_{\Omega_c}]{\sigma\,\partial_t a_z }{a_z' } 
        + \ivol[_\Omega]{\bh }{\boldsymbol{1}_z \times \Grad[]{a_z' } }  
        \\
        = 
        \ivol[_{\Omega_s}]{j_s}{a_z'}, 
        \label{eq:mqs_b_conform_magdyn_a_2D}
    \end{equation}
    for all $a_z' \in \Hone[_0]{\Omega}$. The space $ H^{-1}(\Omega)$ is the dual of $\Hone[_0]{\Omega}$.
\label{prob:multiscale-2D}
\end{prob}
%
%

\subsection{Magnetic flux density conforming formulations: static case}
\label{subsection:mqs_b_static}

The static case can be derived as a particular case of the dynamic problem where eddy currents are neglected. 
The following three-dimensional weak form is obtained from \eqref{eq:mqs_b_conform_magdyn_a}: 
find $\vec{a} \in \Hcurl[_0]{\Omega}$ such that
\begin{gather}
    \ivol[_\Omega]{\bh }{\Curl{\vec{a}'}} 
    = 
    \ivol[_{\Omega_s}]{\vec{j}_s}{\vec{a}'},
    \label{eq:mqs_b_conform_magsta_a}
\end{gather}
for all $\vec{a}' \in \Hcurl[_0]{\Omega}$.
The vector potential $\ba$ is not uniquely defined and a gauge condition must be 
imposed.

Analogously the following two-dimensional weak form is derived from \eqref{eq:mqs_b_conform_magdyn_a_2D}: 
find $a_z \in H_0^1(\Omega)$ such that
\begin{gather}
    \ivol[_\Omega]{\bh }{\boldsymbol{1}_z \times \Grad[]{a_z' } } 
    = 
    \ivol[_{\Omega_s}]{j_s}{a_z'},
    \label{eq:mqs_b_conform_magsta_a_2D}
\end{gather}
for all $a_z' \in \Hone[_0]{\Omega}$

\vspace{10mm}
%
%
\section{Multiscale magnetic induction conforming formulations}
\label{section:models}

A first approach in numerical homogenization consists in precomputing the material 
law. In the case of a material with a linear law and periodic microstructure, only 
one mesoscale problem must be solved in order to get the homogenized quantity 
independent of the macroscale mesh. For the homogenized magnetoquasistatic Problem 
\ref{prob:P_P_homogenized_problem}, the macroscale problem is governed by 
(\ref{eq:equations_P_P_homogenized_macro} a-b)--(\ref{eq:materiallaw_P_P_homogenized_macro} a-b). 
The homogenized magnetic constitutive law (\ref{eq:materiallaw_P_P_homogenized_macro} b)
can be computed by solving the boundary value mesoscale problem \eqref{eq:b_conform_macro_materiallaw_h}. 
For the reversible nonlinear magnetic 
material laws, the points of the material law $\boldsymbol{\mathcal{H}}_M$ can thus be 
computed for different values $\bb_M$, e.g., on the grid 
$\bb_M = (i \, \Delta b_M, k \, \Delta b_M, j \, \Delta b_M)$, with the discretization 
$i, j, k = -N, -(N-1), ... -1, 0, 1, ..., N-1, N$ in each direction and $\Delta b_M = b_M/(2N)$ 
the discretization step and then interpolate to get the values of at any point of 
the application range. This approach was used in \cite{bottauscio-chm-08}.

Hereafter, we develop a coarse-to-fine method inspired by the HMM method
introduced by Weinan E and Enqguist~\cite{Abdulle:2009p1255,Abdulle:2003p808,E:2003p1296,e-hmm-11,E:2003p1295,E:2003p1294,E:2003p1288}. Note
that the so-called FE$^2$ method~\cite{geers-fe2-01,kouznetsova-fe2-01} popular in the
computational mechanics community predates the HMM method and is based on the same
overall philosophy, albeit in a more restrictive setting. This method allows to 
upscale on-the-fly a homogenized material law from the mesoscale problems that 
account for eddy currents at the mesoscale level. These mesoscale problems also 
allow to recover exact electromagnetic fields at the mesoscale level. This approach 
becomes quasi-unavoidable when dealing with problems with hysteresis for which 
the pre-computation of the homogenized magnetic laws described above and the computation 
of local fields are not adapted as they do not account for the history of the material.

In this section we derive the magnetic vector potential formulations for the homogenized 
problem starting with the mesoscale problems governed by the distributional equations 
(\ref{eq:equations_P_P_homogenized_meso} a-b)--(\ref{eq:materiallaw_P_P_homogenized_meso} a-b) 
and the macroscale problem governed by the distributional equations 
(\ref{eq:equations_P_P_homogenized_macro} a-b)--(\ref{eq:materiallaw_P_P_homogenized_macro} a-b).
The index $_m$ is used to denote the restriction of first order terms indexed $_0$ 
on the mesoscale domain $\Omega_m$ (e.g., the restriction of the field $\bb_0$ 
on $\Omega_m$ is denoted by $\bb_m$). 
\deleted{
We also assume the spatial coordinate $\bx$ and the time instance $t$ as parameters 
in mesoscale problems. 
}
The index $M$ refers to the macroscale problems.

\subsection[Magnetic flux density conforming formulations]
{Magnetic flux density conforming multiscale formulations: dynamic case}
\label{sec:phd4_sec2MD_a-v}
\subsubsection{The macroscale problem}
\label{macro-equations}

The macroscale magnetoquasistatic problem was derived in equations
(\ref{eq:equations_P_P_homogenized_macro} a-b) of the homogenized Problem \ref{prob:P_P_homogenized_problem}
%
%
\begin{subequations}
    \begin{equation}
        \Curl[_x]{\bh_M}= \bj_{M}, \quad \Curl[_x]{\be_M}= -\partial_t \bb_M, \quad \Div[_x]{\bb_M}= 0 
        \quad \mathrm{in} \, \,  \Omega \times \mathcal{I}, 
    \tag{\theequation\,a-c}
    \end{equation}
\label{eq:b_conform_macro_maxwell_1}
\end{subequations}
\vspace{-8mm}
\begin{subequations}
    \begin{multline}
        \bh_M(\bx, t) = \boldsymbol{\mathcal{H}}_M \left( \bb_M + \bB_c, \bx \right), 
        \\
        \bj_{M}(\bx, t) = \boldsymbol{\mathcal{J}}_M \left(\be_M + \bE_c, \bx \right)
        \quad \forall (\bx, t) \in \Omega \times \mathcal{I}.	 
    \tag{\theequation\,a-b}
    \end{multline}
\label{eq:b_conform_macro_materiallaw_1}	
\end{subequations}
\vspace{-\baselineskip} 

In (\ref{eq:b_conform_macro_materiallaw_1} a) we use the mapping $\boldsymbol{\mathcal{H}}_M$
instead of the mapping $\boldsymbol{\mathcal{B}}_M$ originally used in Problem \ref{prob:P_P_homogenized_problem}.
This mapping is guaranteed to be uniquely defined if $\boldsymbol{\mathcal{B}}$
is a maximal monotone mapping \cite{visintin-tsh-07-b}.
\noindent The unknown homogenized fields $\vec{h}_M, \vec{b}_M,
\vec{e}_M$ and $\vec{j}_M$ exhibit slow fluctuations; they can therefore be well
approximated on a coarse mesh. The macroscale fields satisfy the same boundary 
conditions as the multiscale fields. Appropriate initial conditions must also be 
provided as specified in \textsc{Assumption} \ref{ass:PH_multiscale_Convergence_IC}. 
Note however that the constitutive laws (\ref{eq:b_conform_macro_materiallaw_1} a-b) 
are not readily available at the macroscale level. They will be upscaled using the 
mesoscale fields.

In the case of a linear electric law $\bj_{M}=\boldsymbol{\mathcal{J}}_M(\be_M) 
= \sigma_M \vec{e}_M$, one computation suffices to extract the homogenized 
conductivity $\sigma_M$ (see details in \cite{bensoussan-ahm-11, niyonzima-chm-13, niyonzima-chm-14}). 
In the case of a nonlinear mapping $\bH_M$, we derive another mesoscale problem 
which accounts for the eddy current effects at the mesoscale level. 
This mesoscale problem (with eddy currents) is thus embedded in a HMM approach
to compute the constitutive homogenized magnetic law on the fly.
Furthermore, it enables the accurate computation of local mesoscale fields
and the upscaling of more accurate global quantities such as the eddy
currents losses. 
\added{
The derivation of the homogenized constitutive laws from the solution of the mesoscale 
time-dependent problem (\ref{eq:equations_P_P_homogenized_meso} a-b)--(\ref{eq:materiallaw_P_P_homogenized_meso} a-b)
instead of the boundary value problem defined by \ref{eq:b_conform_macro_materiallaw_h}
was proved by Visintin (see e.g., \cite[Theorem 7.3]{visintin-tsh-08}).
}

Using results of Section \ref{subsection:mqs_b_dynamic} we can derive the 
three-dimensional macroscale weak formulation of 
\eqref{eq:b_conform_macro_maxwell_1}--\eqref{eq:b_conform_macro_materiallaw_1}.
\begin{prob}[\textbf{Weak form of the three-dimensional MQS macroscale problem}] The weak form of 
the three-dimensional macroscale problem reads:
    find $\ba_M \in \Ltwo{0, T; \Hcurl[_0]{\Omega}}$ with $\partial_t \ba \in \Ltwo{0, T; (\Hcurl[_0]{\Omega})^{\ast}}$ such that
    \vspace{-.5\baselineskip}
    \begin{equation}
    \Big( \sigma_M \partial_t \ba_M , \ba_M' \Big)_{\Omega_{c}} + \Big( \bh_M,  \Curl[_x]{\ba_M' } \Big)_{\Omega}
    = 
    \Big( \bj_s , \ba_M' \Big)_{\Omega_s},
    \label{eq:mqs_b_conform_magdyn_a_multiscale}
    \end{equation}
    hold for all test functions $\ba_M' \in \Hcurl[_0]{\Omega}$. 
    \label{eq:mqs_b_conform_magdyn_a_macroscale_3D}
\end{prob}
The macroscale magnetic field $\bh_M(\bx, t) = \boldsymbol{\mathcal{H}}_M(\Curl[_x]{\ba_M} + \bb_c, \bx, t)$
is dependent on the mesoscale solutions $\bb_{\mathrm{c}}$.
The vector 
\begin{equation}
  \bb_{\mathrm{c}} = (\bb_{\mathrm{c}}^{(1)},\; \bb_{\mathrm{c}}^{(2)},\; \ldots,\; \bb_{\mathrm{c}}^{(i)},\; \ldots,\; \bb_{\mathrm{c}}^{(\text{N}_{\mathrm{GP}})})
\end{equation}
is a collection of magnetic field corrections obtained by applying the solution 
operators in \eqref{eq:solution-operators-b} for mesoscale problems corresponding 
to Gau{\ss} points $\bx^{(i)}$. The vector $\vec{j}_M$ represents the eddy 
currents and $\vec{j}_s$ represents 
the source current density imposed in the inductors $\Omega_s$. 
\deleted{The macroscale 
domain $\Omega$ (resp. $\Omega_c$) can be divided into the multiscale domain 
$\Omega^h$ (resp. $\Omega_c^h$) where the homogenization is performed and a 
non-multiscale domain $\Omega^{nh}$ (resp. $\Omega_c^{nh}$) where classical weak 
formulations hold.
}

\noindent For the two-dimensional case, we get the following problem:
\begin{prob}[\textbf{Weak form of the two-dimensional MQS macroscale problem}] The weak form of 
    the two-dimensional macroscale problem reads:
    find $a_{zM} \in \Ltwo{0, T; \Hone[_0]{\Omega}}$ with $\partial_t a_{zM} \in \Ltwo{0, T; H^{-1}(\Omega)}$ such that
    \vspace{-.5\baselineskip}
    \begin{equation}
        \Big( \sigma_M \partial_t a_{zM} , a_{zM}' \Big)_{\Omega_{c}} + 
        \Big( \bh_M, \boldsymbol{1}_z \times \Grad[_x]{a_{zM}' }  \Big)_{\Omega} 
        =
        \Big( j_s , a_{zM}' \Big)_{\Omega_s},
    \label{eq:mqs_b_conform_magdyn_a_multiscale_2D}
    \end{equation}
    hold for all test functions $a_{zM}' \in H_0^1(\Omega)$.
    \label{eq:mqs_b_conform_magdyn_a_macroscale_2D}
\end{prob}

The homogenized magnetic law $\boldsymbol{\mathcal{H}}_M$ in equations
\eqref{eq:mqs_b_conform_magdyn_a_multiscale} for the three-dimensional problem and in 
\eqref{eq:mqs_b_conform_magdyn_a_multiscale_2D} for the two-dimensional problem is 
upscaled using the mesoscale fields as described in the following section.

\subsubsection{The mesoscale problem}
\label{sec:phd4_b-form_meso}

The governing equations of the mesoscale problem with eddy currents which, unlike 
problem \eqref{eq:b_conform_macro_materiallaw_h}--\eqref{eq:b_conform_macro_materiallaw_j}, 
also enables to recover accurate local electromagnetic fields, are a modified version 
of the two-scale problem (\ref{eq:equations_P_P_homogenized_meso} a-b)--(\ref{eq:materiallaw_P_P_homogenized_meso} a-b). 
These equations read:
\begin{subequations}
    \begin{equation}
        \Curl[]{} \bh_m^{\varepsilon}= \bj_{m}, 
        \qquad \Curl[]{_x} \be_M+ \Curl[]{_y} \be_1 = -\partial_t\bb_m, 
    \tag{\theequation\,a-b}
    \end{equation}
\label{eq:b_conform_meso_maxwell_1}
\end{subequations}
\vspace{-8mm}
\begin{subequations}
    \begin{equation}
        \bh_m(\bx, \by, t) = \boldsymbol{\mathcal{H}}(\bb_m(\bx, \by, t), \bx, \by),
        \quad 
        \bj_{m}(\bx, \by, t) = \boldsymbol{\mathcal{J}}(\be_m(\bx, \by, t), \bx, \by), 
    \tag{\theequation\,a-b}
    \end{equation}
\label{eq:b_conform_meso_materiallaw_1}	
\end{subequations}

\vspace{-\baselineskip}
\noindent in which we keep the $\CurlSymb$ of $\bh^{\varepsilon}$ instead of using 
its two-scale decomposition given in (\ref{eq:equations_P_P_homogenized_meso} a). In this
equation, $\vec{h}_m^{\varepsilon}$ is the restriction of the multiscale
magnetic field $\vec{h}^{\varepsilon}$ to the representative volume element
$\Omega_m$ also called ``mesoscale domain''. We can thus use both
nonlinear reversible and irreversible (hysteretic) material
laws. Problem (\ref{eq:b_conform_meso_maxwell_1} a-b)--(\ref{eq:b_conform_meso_materiallaw_1} a-b)
contains macroscale fields assumed constant at the mesoscale
level, so that the mesoscale problem can be written in terms of the
mesoscale coordinates $\by$. This is the case if the scale separation assumption 
is fulfilled.

The two-scale convergence theory allows us to express the $\CurlSymb$ of the
electric field at the mesoscale level in terms of the $\CurlSymb$ of the electric
field at the macroscale and the $\CurlSymb$ of the mesoscale correction
term i.e.\ $\Curl[_y]{\vec{e}_m} = \Curl[_x]{\vec{e}_M}+
\Curl[_y]{\vec{e}_1}$. Using the Faraday law at the macroscale together with
the vector identity $\nolinebreak{\Curl[_y]{(\partial_t b_M \times \by)} = 
(n-1) \partial_t b_M}$ ($n = 2, 3$ for two-dimensional and three-dimensional problems, respectively)
we can write:
\begin{multline}
    \Curl[_y]{\vec{e}_m} = \Curl[_y]{\Big(\be_1 + \be_M + \kappa (\Curl[_y]{\vec{e}_M} \times \vec{y}) \Big)} 
    \\
    = \Curl[_y]{\Big(\be_1 + \be_M - \kappa (\partial_t \bb_M \times \by) \Big)}
\label{eq:electricField_meso}
\end{multline}
with $\kappa = (n-1)^{-1}$, since $\Curl[_y]{\vec{e}_M}  \equiv 0$. Similar
developments have been proposed in \cite{meunier-chm-10} and \cite{feddi-thm-97} 
for the electric and the magnetic fields in linear cases. Inserting the orthogonal 
decomposition of the mesoscale magnetic induction $\bb_m = \bb_M +\Curl[_y]{\ba_c}$ 
we get
\begin{equation}
    \Curl[]{_x} \be_M+ \Curl[]{_y} \be_1 = -\partial_t(\bb_M + \Curl[_y]{\ba_c} ).
\label{eq:electricField_meso_1}
\end{equation}
From \eqref{eq:electricField_meso_1} we get $\Curl[_y]{(\be_1 + \partial_t \ba_c)} = \boldsymbol{0}$
which, together with the orthogonal decomposition \eqref{eq:space-decomposition_Hdiv} 
leads to the expression of the first order term of the electric field $\be_1$ in 
terms of the correction terms $\ba_c$ and $v_c$ as:
\begin{equation}
    \be_1 = -\partial_t \ba_c - \Grad[_y]{v_c}.
\label{eq:electricField_e1}
\end{equation}
At the mesoscale level, the first order term $\vec{e}_1(\bx, \cdot, t)$ must be chosen in
$\Hcurl[_*]{ \mathcal{Y}}$ for almost every $(\bx, t) \in \mathbb{R}_T^3$. 
In Section \ref{scale_trans} we will show that $\vec{a}_c$ is tangentially periodic 
and we will choose $v_c$ to be periodic on the mesoscale domain $\Omega_m$.
Using these developments, we can derive the mesoscale three-dimensional weak formulation.
\begin{prob}[\textbf{Weak form of the three-dimensional MQS mesoscale problem}] The weak form of 
    the three-dimensional mesoscale problem reads: 
    find $\ba_c \in \Ltwo{0, T; \Hcurl[_*]{\mathcal{Y} } }$ with 
    $\partial_t \ba_c \in \Ltwo{0, T; (\Hcurl[_*]{\mathcal{Y} })^{\ast} }$ and 
    $v_c \in \Ltwo{0, T; \Hone[_*]{\mathcal{Y} } }$ such that
    \begin{multline}
        \Big( \sigma \partial_t \ba_c , \ba_c' \Big)_{\Omega_{mc}  } +
        \Big( \bh ,  \CurlSymb_y \ba_c' \Big)_{\Omega_{m} } + 
        \\
        \Big( \sigma \GradSymb_y v_c , \ba_c' \Big)_{\Omega_{mc} } = 
        \Big( \sigma (\be_M -\kappa \partial_t \bb_M \times \by ), \ba_c' \Big)_{\Omega_{mc} },
        \label{eq:weakform_magdyn_vecPot_micro}
        \end{multline}
        \begin{multline}
        \Big( \sigma \partial_t \ba_c , \GradSymb_y v_c' \Big)_{\Omega_{mc} } + 
        \Big( \sigma \GradSymb_y v_c ,  \GradSymb_y v_c' \Big)_{\Omega_{mc} } =
        \\
        \Big( \sigma (\be_M - \kappa \partial_t \bb_M \times \by), \GradSymb_y v_c'
        \Big)_{\Omega_{mc} } + \Big< \bn\cdot \bj_{M}, v_c' \Big>_{\Gamma_{gm} } 
    \label{eq:weakform_magdyn_scalarPot_micro}
    \end{multline}
    hold for all test functions $\ba_c' \in \Hcurl[_*]{\mathcal{Y} }$ and
    $v_c' \in \Hone[_*]{\mathcal{Y} }$.
    \label{eq:mqs_b_conform_magdyn_a_mesoscale_3D}
\end{prob}
The magnetic field is given by $\bh(\bx, \by, t) = \boldsymbol{\mathcal{H}}(\CurlSymb_y \ba_c(\bx, \by, t) + \bb_M(\bx, t), \bx, \by)$.
and the boundary term $\left< \bn \times \bh, \ba_c^{'} \right>_{\Gamma_m}$ is omitted 
due to the periodicity of $\bh = \bh_0$ (see the definition of function space 
in \eqref{eq:eq_Meso_Space}) and of $\ba_c^{'}$. The domain $\Omega_{mc}$ with 
boundary $\Gamma_{gm}$ is the conducting part of the mesoscale domain and the electric 
current density $\bj_{M} = \sigma_M \be_M$ is obtained from the macroscale solution.

For the two-dimensional case, the following mesoscale weak formulation can be derived.
\begin{prob}[\textbf{Weak form of the two-dimensional MQS mesoscale problem}] The weak form of 
the two-dimensional mesoscale problem reads: 
    find $a_{zc} \in \Ltwo{0, T; \Hone[_*]{\mathcal{Y} } }$ with 
    $\partial_t a_{zc} \in \Ltwo{0, T; (\Hone[_*]{\mathcal{Y} })^{\ast} }$
    and $u_c$ piecewise constant on $\Omega_{mc}$ for almost every $(\bx, t) \in \mathbb{R}_T^3$ such that
    \begin{multline}
        \Big( \sigma \partial_t a_{zc} , a_{zc}' \Big)_{\Omega_{mc}  } +
        \Big( \bh, \boldsymbol{1}_z \times \Grad[_y]{a_{zc}' } \Big)_{\Omega_{m} } + 
        \\
        \Big( \sigma u_c, a_{zc}' \Big)_{\Omega_{mc} } = 
        \Big( \sigma (\be_M -\kappa \partial_t \bb_M \times \by ), \boldsymbol{1}_z a_{zc}' \Big)_{\Omega_{mc} },
    \label{eq:weakform_magdyn_vecPot_micro_2D}
    \end{multline}
    \begin{multline}
        \Big( \sigma \partial_t a_{zc}, u_c' \Big)_{\Omega_{mc} } \! + 
        \Big( \sigma u_c,  u_c' \Big)_{\Omega_{mc}} \!
        = \Big( \sigma (\be_M - \kappa \partial_t \bb_M \times \by), \boldsymbol{1}_z
        u_c' \Big)_{\Omega_{mc}} \label{eq:weakform_magdyn_scalarPot_micro_2D}
    \end{multline}
    hold for all test functions $a_{zc}' \in \Hone[_*]{\mathcal{Y} }$ and $u_c'$
    piecewise constant on $\Omega_{mc}$.
    \label{eq:mqs_b_conform_magdyn_a_mesoscale_2D}
\end{prob}
%
%

\subsubsection{Scale transitions}
\label{scale_trans}
The macroscale and the mesoscale problems in Sections \ref{macro-equations}
and \ref{sec:phd4_b-form_meso} are not yet well-defined. Indeed, the macroscale
magnetic law $\boldsymbol{\mathcal{H}}_M$ is not readily available at the
macroscale level and the mesoscale problem requires source terms $\vec{b}_M,
\vec{e}_M$ and $\vec{j}_M$ and proper boundary conditions to be
well-posed. 
These two problems need to fill the missing information by exchanging
data between the macro and meso levels. The so-called scale transitions
comprise the \textit{downscaling} and the \textit{upscaling} stages (see Figure~\ref{FE2_ppe}).

During the \textit{downscaling}, the macroscale fields are imposed as source
terms for the mesoscale problem.  Boundary conditions for the mesoscale
problem are also determined so as to respect the two-scale convergence of
the physical fields, i.e., the convergence of the magnetic flux density $\langle
\bb_m \rangle_{\Omega_m} = \bb_M$ leads to the following condition on the
tangential component of the correction term of the magnetic vector potential
$\ba_c$, which is fulfilled if
\begin{equation}
    \int_{\Omega_m} \Curl[]{\ba_c(\bx, \by, t)} \textrm{d} y   = 
    \oint_{\Gamma_m} \bn \times \ba_c(\bx, \by, t)\textrm{d} y = \boldsymbol{0}, \qquad \forall (\bx, t) \in \mathbb{R}_T^3.
    \label{eq:conv-magfluxdens}
\end{equation}
This condition is fulfilled if $\vec{a}_c (\bx, \cdot, t)$ belongs to the space
$\Hcurl[_*]{\mathcal{Y} }$, i.e. if $\vec{a}_c$ is tangentially periodic on the cell.
\begin{figure}[t]
  \begin{center}
    \includegraphics[width=0.6\textwidth]{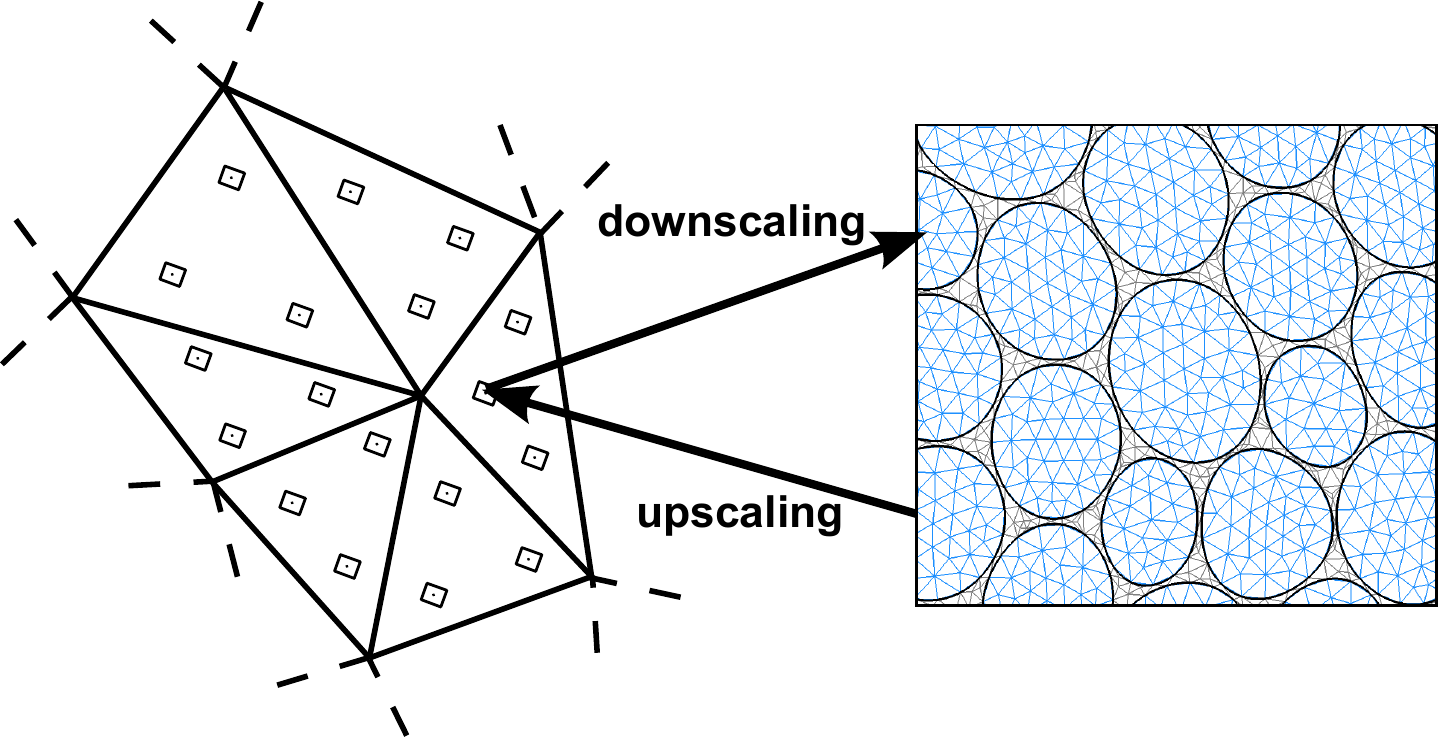}
  \end{center}
  \caption{Scale transitions between the macroscale (left) and the mesoscale (right)
    problems. Downscaling (macro to meso): obtaining proper boundary conditions
    and the source terms for the mesoscale problem from the current macroscale
    solution. Upscaling (meso to macro): calculating effective quantities (e.g.\
    material properties) for the macroscale problem from the mesoscale
    solution}~\cite{niyonzima-chm-12}.
  \label{FE2_ppe}
\end{figure} 
Additionally, $\Grad[_y]{v_c}(\bx, \cdot, t) = \vec{e}_1(\bx, \cdot, t) - \partial_t \vec{a}_c(\bx, \cdot, t)$ also
belongs to $\Hcurl[_*]{\mathcal{Y} }$, which is automatically ensured by the \textit{curl theorem}:
\begin{equation}
    \int_{\Gamma_m } \vec{n} \times \Grad[_y]{v_c} \text{d} y = \int_{\Omega_m } \Curl[_y]{\Grad[_y]{v_c}} \text{d} y.
    \label{eq:periodicityOf_vc}
\end{equation} 
Further we choose a periodic $v_c$.

The convergence of the electric current density $\langle \bj_m \rangle_{\Omega_{mc}} = \bj_M$
also leads to the following relation:
\begin{multline}
    \int_{\Omega_{mc}} \bj_c(\bx, \by, t) \, \textrm{d} y 
    \\
    = - \int_{\Omega_{mc}} \sigma \Big(\partial_t \ba_c(\bx, \by, t) +
    \Grad[]{v_c(\bx, \by, t)}\Big) \, \textrm{d} y =  \boldsymbol{0}, \,\, \forall (\bx, t) \in \mathbb{R}_T^3.
    \label{eq:conv-elec_flux_field}
\end{multline}

The \emph{upscaling} consists in computing the missing constitutive laws 
$\sigma_M$, $\displaystyle \boldsymbol{\mathcal{H}}_M$ together with 
$\partial\boldsymbol{\mathcal{H}}_M/\partial \bb_M$ at the macroscale using the 
mesoscale fields. Due to the linearity of the electric law, the asymptotic expansion
theory can be applied. Therefore, we compute once and for all the homogenized electric
conductivity by solving a unique cell problem. A similar approach was 
also adopted in \cite{bottauscio-chm-08}.

The upscaling of the nonlinear magnetic law is performed by
averaging the magnetic field (consequence of the two-scale convergence of
the magnetic field):
\begin{multline}
\bh_M(\bx, t) = \bH_M\left( \bb_M(\bx, t)\right) 
\\
= \dfrac{1}{|\Omega_{mc}|} \int_{\Omega_{mc}} \bH\left(\Curl[_x]{\ba_M(\bx, t)} + \Curl[_y]{\ba_c(\bx,
\by, t)}, \bx\right) \textrm{d}y.
\end{multline}
For the $i$-th Gau{\ss} point, the Jacobian expression reads
\begin{equation}
  \frac{\mathrm{d} \bH_{\mathrm{M}}}{\mathrm{d} \bb_{\mathrm{M}}}
  =
  \frac{1}{|\Omega_{\textrm{m}}|} \int_{\Omega_{\textrm{m}}} 
  \left( 
    \dfrac{\partial \bH}{ \partial \bb_{\mathrm{M}}} 
    +
    \dfrac{\partial \bH}{ \partial \bb_{\mathrm{c}}}   
    \dfrac{\partial \bB_{\mathrm{c}}}{ \partial \bb_{\mathrm{M}}}
  \right) \text{d} \by
\label{eq:Homogenized-Law-FE-HMM-1}
\end{equation}
with $\bb_{\mathrm{c}}(\bx, \by, t) = \Curl[_y]{\ba_{\mathrm{c}}(\bx, \by, t)} 
= \bB_{\mathrm{c}}\left(\by, \Curl[_x]{\ba_{\mathrm{M}}}(\bx, t)\right)$ 
and $\bb_M = \Curl[_x]{\ba_{\mathrm{M}}}$. The derivative w.r.t. the mesoscale vector potential $\ba_{\mathrm{M}}$ 
is given by
\begin{equation}
  \dfrac{\mathrm{d} \bH_{\mathrm{M}}}{\mathrm{d} \ba_{\mathrm{M}}}
  =
  \dfrac{\mathrm{d} \bH_{\mathrm{M}}}{\mathrm{d} \bb_{\mathrm{M}}}
  \dfrac{\mathrm{d} \bb_{\mathrm{M}}}{\mathrm{d} \ba_{\mathrm{M}}}.
\end{equation}

The computation \eqref{eq:Homogenized-Law-FE-HMM-1} involves the Fr\'{e}chet 
derivative of $\bB_{\mathrm{c}}$ with respect to the macroscale magnetic density 
$\bb_{\mathrm{M}}$. This derivative can be evaluated numerically using the finite 
difference. In \cite{niyonzima-chm-14}, several mesoscale problems  
per Gau{\ss} point were solved in parallel. A first problem is solved using 
\eqref{eq:weakform_magdyn_vecPot_micro}--\eqref{eq:weakform_magdyn_scalarPot_micro}
for the three-dimensional problems (resp. 
\eqref{eq:weakform_magdyn_vecPot_micro_2D}--\eqref{eq:weakform_magdyn_scalarPot_micro_2D}
for the two-dimensional case) to find the solution when a macroscale source $\vec{b}_M$ is applied. 
Then, a time and space independent magnetic induction perturbation
term $\delta \bb_i$ oriented along the $i$ directions ($_i = _x, _y$ and $_z$) 
is added to the macroscale source terms. Therefore, three (resp. two) additional problems analogous to    
\eqref{eq:weakform_magdyn_vecPot_micro}--\eqref{eq:weakform_magdyn_scalarPot_micro}
(resp. \eqref{eq:weakform_magdyn_vecPot_micro_2D}--\eqref{eq:weakform_magdyn_scalarPot_micro_2D}
for the two-dimensional case) are solved in order to determine the Jacobian 
$\mathrm{d} \bH_{\mathrm{M}}/\mathrm{d} \bb_{\mathrm{M}}$ needed for the 
Newton-Raphson scheme. 
The total magnetic induction $\bb_m$ for these problems are expressed as:
\begin{equation}
\bb_m = \vec{b}_M + \Curl[_y]{\vec{a}_c} + \delta \bb_i = 
\Curl[_y]{\left( \vec{a}_c + \kappa(\vec{b}_M \times \by) + \kappa(\delta \bb_i \times \by) \right)}
 = \Curl[_y]{\ba_m},
\end{equation}
which can be derived from the total magnetic vector potential: 
\begin{equation}
\ba_m = \vec{a}_c - \Grad[_y]{v_c} + \kappa(\vec{b}_M \times \by) + \kappa(\delta \bb_i \times \by).
\label{eq:weakform_magdyn_vecPot_micro_pert_a}
\end{equation}
These developments allow to transform the three dimensional equation 
\eqref{eq:weakform_magdyn_vecPot_micro} into
\begin{multline}
\Big( \sigma \partial_t \ba_c , \ba_c' \Big)_{\Omega_{mc}  } +
\Big( \boldsymbol{\mathcal{H}}(\CurlSymb_y \ba_c + \bb_M + \delta \bb_i, \bx) ,
\CurlSymb_y \ba_c' \Big)_{\Omega_{m} } + 
\\
\Big( \sigma \GradSymb_y v_c , \ba_c' \Big)_{\Omega_{mc} } = 
\Big( \sigma (\be_M -\kappa \partial_t \bb_M \times \by ), \ba_c' \Big)_{\Omega_{mc} }.
\label{eq:weakform_magdyn_vecPot_micro_pert}
\end{multline}
Notice that the time derivative of the constant term in equation
\eqref{eq:weakform_magdyn_vecPot_micro_pert_a} disappears. We also
modify the two dimensional equation \eqref{eq:weakform_magdyn_vecPot_micro_2D} as 
\begin{multline}
\Big( \sigma \partial_t a_{zc} , a_{zc}' \Big)_{\Omega_{mc}  } +
\Big( \boldsymbol{\mathcal{H}}(\boldsymbol{1}_z \times \Grad[_y]{a_{zc}} +
\bb_M + \delta \bb_i, \bx) , \boldsymbol{1}_z \times \Grad[_y]{a_{zc}' }
\Big)_{\Omega_{m} } + 
 \\
\Big( \sigma u_c, a_{zc}' \Big)_{\Omega_{mc} } = 
\Big( \sigma (\be_M -\kappa \partial_t \bb_M \times \by ), \boldsymbol{1}_z a_{zc}' \Big)_{\Omega_{mc} }.
\label{eq:weakform_magdyn_vecPot_micro_2D_333}
\end{multline}
Equations \eqref{eq:weakform_magdyn_scalarPot_micro} and \eqref{eq:weakform_magdyn_scalarPot_micro_2D} remain unchanged.
This leads to the solution $\vec{h}_M + \delta_{\bb_i } \vec{h}_M$ where
$\delta_{\bb_i } \vec{h}_M$ is the perturbation of the magnetic field in
direction $i$. We can therefore compute the elements of the tangent matrix as:
\begin{equation}
\left(\dfrac{\partial \boldsymbol{\mathcal{H}}_M }{\partial \vec{b}_M}
\right)_{ij} \approx 
\dfrac{(\delta_{\bb_{i} } \vec{h}_M)_j}{\delta \bb_{i}}.
\end{equation}
Further mathematical justifications of the numerical computation of the tangent matrix 
are given in Section \ref{fe-hmm-discrete}.

\subsection{Finite element implementation}
\label{fe-hmm-discrete}
\added{In this section we discuss the numerical implementation of the 
homogenized problem using the finite element method. The  numerical approximation
involves errors the sources of which can be numerous in the case of the MQS problem:
}
\begin{itemize}
    \item the error due to the the finiteness of the mesoscale domain (instead of $\varepsilon \rightarrow 0$), 
    \item the error due to the modified mesoscale problem (\ref{eq:b_conform_meso_maxwell_1} a-b),
    \item the error due to scale transitions,
    \item the error due to the approximation using a finite dimensional space in the Galerkin approximation,
    \item the error due to Euler's time stepper,
    \item the error in the Newton--Raphson scheme used for solving the nonlinear macroscale and the mesoscale problems,
    \item the error due to the resolution of the linear systems, 
    \item the error in reduced Jacobian,
    \item the error resulting from the application of the homogenization near the boundaries of the computational domain, etc.
\end{itemize}
This paper does not deal with the error analysis.
\deleted{For the sake of simplicity, we consider a modified magnetic 
vector potential so that $\be=-\partial_t \ba$, i.e.\ we do not include an 
electric scalar potential $v$ in the macroscale and mesoscale problems. Note that 
we abuse notation using the same quantity $\ba$. The extension to the $\ba-v$ formulation 
is straightforward.
}

Using a similar approach to the one used in \cite{niyonzima-chm-16}, the macroscale and mesoscale 
equations are solved using the finite element method. The fields $\ba_{\mathrm{M}}^H$ 
and $\ba_{\mathrm{c}}^H$ are approximation of the continuous fields $\ba_{\mathrm{M}}$ and 
$\ba_{\mathrm{c}}$ on the discretized computational domain and $\ba_{\mathrm{c}}^H \in (0, T ] \times \bW_{H,0}^\mathrm{M}$ 
and $\ba_{\mathrm{c}}^H \in (0, T ] \times \bW_{h,0}^\mathrm{m}$ where $\times \bW_{H,0}^\mathrm{M}$ 
and $\bW_{h,0}^\mathrm{m}$ are discrete subspaces of $\Hcurl[_e^0]{\Omega}$ and $\Hcurl[_e^0]{\Omega_m}$
\begin{multline}
    \ba_{\mathrm{M}}(\bx,t) \approx \ba_{\mathrm{M}}^H(\bx, t) = 
    \sum_{p=1}^{N_\mathrm{M}} \ah_{\textrm{M},p}(t) \ba^{\prime}_{\textrm{M},p}(\bx) 
    \\
    \qquad \textrm{and} \qquad
    \ba_{\mathrm{c}}^{(i)}(\by,t) \approx \ba_{\mathrm{c}}^H(\by, t) = 
    \sum_{p=1}^{N_\mathrm{c}} \ah_{\textrm{c},p}^{(i)}(t) \ba^{\prime}_{\textrm{c},p}(\by),
\label{eq:fem_macro_meso}
\end{multline}
where the superscript $i = 1, 2, \cdots, N_\textrm{GP}$ refers to the enumeration 
of mesoscale problems, $N_\mathrm{M}$ and $N_\mathrm{c}$ are the number of degrees 
of freedom for discretized fields at the macroscale and the mesoscale, respectively.
Space discretization leads to the semidiscrete coupled problem:

Find waveforms $[\bah_{\mathrm{M}}(t),\bah_{\mathrm{c}}^{(1)}(t), \ldots, \bah_{\mathrm{c}}^{(N_\mathrm{GP})}(t)]$
such that
\begin{equation}
    \bM_{\mathrm{M}}\partial_t \bah_\mathrm{M}+\bF_{\mathrm{M}}(\bah_\mathrm{M},\bah_{\mathrm{c}})=0, 
  \label{eq:newton-macro-semidiscrete}
\end{equation}
and for the mesoscale problems $i=1,\ldots,N_\mathrm{GP}$
\begin{equation}
    \bM_{\mathrm{m}}\partial_t \bah_{\mathrm{c}}^{(i)}+\bF_\mathrm{m}(\bah_{\mathrm{c}}^{(i)},\bah_{\mathrm{M}}^{(i)}, \partial_t\bah_\mathrm{M}^{(i)})=0
    \label{eq:newton-meso-semidiscrete}
\end{equation}
for a given set of initial values $[\bah_{\mathrm{M}}(t_0),\bah_{\mathrm{c}}^{(1)}(t_0), \ldots, \bah_{\mathrm{c}}^{(N_\mathrm{GP})}(t_0)]$.
where $\bM_{\mathrm{M}} :=\Big( \sigma_\mathrm{M} \ba_\mathrm{M}, \ba_\mathrm{M}^{\prime}\Big)_{\Omega_\mathrm{c}}$ with $\ba_\mathrm{M}$ and 
$\ba_\mathrm{M}^{\prime}$, the ansatz functions, 
the functions $\bF_\mathrm{M}(\cdots)$ 
and $\bF_\mathrm{m}(\cdots)$ are the semi-discreet terms involving the nonlinear 
magnetic terms stemming from \eqref{eq:mqs_b_conform_magdyn_a_multiscale} and
\eqref{eq:weakform_magdyn_vecPot_micro} by inserting \eqref{eq:fem_macro_meso}.

The time discretization using an implicit Euler method 
\deleted{leads to the system
(for each time step $k=0,\ldots,N_\mathrm{TS}$):
find a series of solutions
}
followed by the use of the Newton--Raphson method to solve the resulting nonlinear 
problem leads the following Jacobian 
\begin{multline}
J_{\bR}^{(j,k)} := 
\frac{1}{\Delta t_k}
\begin{pmatrix} 
    \displaystyle{\bM_{\mathrm{M}}\vphantom{\frac{\partial \bF^{(j,k)}_{\mathrm{M}}}{\partial \bah_{\mathrm{M}}^{(j,k)}}}}& 
    \displaystyle{0} & 
    \displaystyle{\cdots} & 
    \displaystyle{0} 
    \\ 
    \displaystyle{0} & 
    \displaystyle{\bM_{\mathrm{m}}\vphantom{\frac{\partial \bF^{(j,k)}_{\mathrm{M}}}{\partial \bah_{\mathrm{M}}^{(j,k)}}}} & 
    \displaystyle{0} & 
    \displaystyle{0}    
    \\ 
    \displaystyle{\vdots} & 0 & \displaystyle{\ddots} & 0 
    \\ 
    \displaystyle{0} & 
    0 & 
    0 & 
    \displaystyle{\bM_{\mathrm{m}}\vphantom{\frac{\partial \bF^{(j,k)}_{\mathrm{M}}}{\partial \bah_{\mathrm{M}}^{(j,k)}}}}
\end{pmatrix}%
+
\\
\begin{pmatrix} 
    \displaystyle{\frac{\partial \bF^{(j,k)}_{\mathrm{M}}}{\partial \bah_{\mathrm{M}}^{(j,k)}}}& 
    \displaystyle{\frac{\partial \bF^{(j,k)}_{\mathrm{M}}}{\partial \bah_{\mathrm{c}}^{(1,j,k)}}} & 
    \displaystyle{\cdots} & 
    \displaystyle{\frac{\partial \bF^{(j,k)}_{\mathrm{M}}}{\partial \bah_{\mathrm{c}}^{(N_{\mathrm{GP}},j,k)}}} 
    \\ 
    \displaystyle{\frac{\partial \bF^{(1,j,k)}_{\mathrm{m}}}{\partial \bah_{\mathrm{M}}^{(1, j,k)}}} & 
    \displaystyle{\frac{\partial \bF^{(1,j,k)}_{\mathrm{m}}}{\partial \bah_{\mathrm{c}}^{(1,j,k)}}} & 
    \displaystyle{0} & 
    \displaystyle{0}    
    \\ 
    \displaystyle{\vdots} & 0 & \displaystyle{\ddots} & 0 
    \\ 
    \displaystyle{\frac{\partial \bF^{(N_{\mathrm{GP}},j,k)}_{\mathrm{m}}}{\partial \bah_{\mathrm{M}}^{(N_{\mathrm{GP}},j,k)}}} & 
    0 & 
    0 & 
    \displaystyle{\frac{\partial \bF^{(N_{\mathrm{GP}},j,k)}_{\mathrm{m}}}{\partial \bah_{\mathrm{c}}^{(N_{\mathrm{GP}},j,k)}}}
\end{pmatrix}%
\label{eq:Jacobian-full}
\end{multline}%
where $\bah_{\mathrm{M}}^{(j,k)}$ and $\bah_{\mathrm{c}}^{(i,j,k)}$ denote the $j^{\mathrm{th}}$
Newton--Raphson iterates and 
\begin{multline}
    \bF_{\mathrm{M}}^{(j, k)} := \bF_{\mathrm{M}} \Big(\bah_{\mathrm{M}}^{(j, k)}, 
    \bah_{\mathrm{c}}^{(j, k)} \Big) 
    \\
    \bF_{\mathrm{m}}^{(i, j, k)} := \bF_{\mathrm{m}} \left(\bah_{\mathrm{c}}^{(i, j, k)}, \bah_{\mathrm{M}}^{(i, j, k)},
    \frac{\bah_{\mathrm{M}}^{(i, j, k)} - \bah_{\mathrm{M}}^{(i, j, k-1)}}{\Delta t_k} \right).
\end{multline}
where the superscripts $k$ and $j$ are used for time steps and the Newton--Raphson 
iterations. See \cite[Section 4]{niyonzima-chm-16} for more details on the derivation 
of the Jacobian \eqref{eq:Jacobian-full}. In practice, one does not solve the system 
above but the Schur complement system with the reduced Jacobian defined by
\begin{equation}
    \bar{J}_{\bR}^{(j,k)} := \frac{\bM_{\mathrm{M}}}{\Delta t_k} + 
    \dfrac{\partial \bF^{(j,k)}_{\mathrm{M}}}{\partial \bah_{\mathrm{M}}^{(j,k)}}
    -\sum_{i=1}^{N_{\mathrm{GP}}}
    \left(
    \dfrac{\partial \bF^{(j,k)}_{\mathrm{M}}}{\partial \bah_{\mathrm{c}}^{(i,j,k)}}  
    \Big(\frac{\bM_{\mathrm{m}}}{\Delta t_k} + 
        \dfrac{\partial \bF^{(i,j,k)}_{\mathrm{m}}}{\partial \bah_{\mathrm{c}}^{(i,j,k)}}
    \Big)^{-1} \, 
    \dfrac{\partial \bF^{(i,j,k)}_{\mathrm{m}}}{\partial \bah_{\mathrm{M}}^{(i,j,k)}}
    \right).
\label{eq:Jacobian-schur-complement-macro}
\end{equation}
%
%

%
%
\begin{center}
\begin{algorithm}[!]
\caption{Pseudocode for the monolithic FE-HMM}\label{alg:FE-HMM-Macro}
\begin{algorithmic}
\INPUT macroscale source $\bj_{\mathrm{s}}$ and mesh.
\OUTPUT macroscale fields, mesoscale fields and global quantities.
\Procedure{macroscale problem}{} 
    \State $t \gets t_{0}$, initialize the macroscale field $\ba_{\mathrm{M}}|_{t_0} = \ba_{\mathrm{M}0}$,
    \For{$(k \gets 1$ To $N_{\mathrm{TS}} )$}              \Comment{\emph{the macroscale time loop $($index $k$$)$}}
    \For{$(j \gets 1$ To $N_{\mathrm{NR}}^{\mathrm{M}} )$} \Comment{\emph{the macroscale NR loop $($index $j$$)$}}
    \For{$(i \gets 1$ To $N_{\mathrm{GP}} )$}              \Comment{\emph{parallel solutions of meso-problems $($index $i$$)$}}
    \State downscale the macroscale sources,
    \State compute the mesoscale fields, see Algorithm~\ref{alg:FE-HMM-Meso},
    \State compute and upscale the homogenized law $\bH_{\mathrm{M}}$ and 
    \State $\partial \boldsymbol{\mathcal{H}}_{\mathrm{M}}/ \partial \bb_{\mathrm{M}}$
    \EndFor
    \State assemble the Jacobian $\bar{J}_{\bR}^{(j,k)}$ from \eqref{eq:Jacobian-schur-complement-macro} to solve the macroscale problem,
    \EndFor
    \EndFor
\EndProcedure
\end{algorithmic}
\end{algorithm}
\end{center}
\begin{center}
\begin{algorithm}[!]
\caption{Pseudocode for one mesoscale problem}\label{alg:FE-HMM-Meso}
\begin{algorithmic}
\INPUT macroscale sources and the mesoscale mesh.
\OUTPUT homogenized law $\bH_{\mathrm{M}}$ and $\partial \boldsymbol{\mathcal{H}}_{\mathrm{M}}/ \partial \bb_{\mathrm{M}}$, per Gau\ss \, point for $N_m$ problems.
\Procedure{mesoscale problem}{} 
    \State prescribe periodic boundary conditions, impose sources,
    \State $t \gets t_{\mathrm{M}}$, initialize the correction term $\ba_{\mathrm{c}}|_{t_{\mathrm{M}}}$
    \For{$(p \gets 1$ To $N^\mathrm{m}_\mathrm{dim})$}  \Comment{\emph{solve $N^\mathrm{m}_\mathrm{dim}$ mesoscale problems for the $k^{th}$ time step} }
    \For{$(j \gets 1$ To $ N_{\mathrm{NR}}^{\mathrm{m}} )$} \Comment{\emph{the mesoscale NR loop $($index $j$$)$}}
    \State assemble the matrix and solve the mesoscale problem.
    \EndFor    
    \EndFor
\EndProcedure
\end{algorithmic}
\end{algorithm}
\end{center}
The overall FE-HMM method is described in Algorithm \ref{alg:FE-HMM-Macro} 
and Algorithm \ref{alg:FE-HMM-Meso}. It starts with the 
initialization of the macroscale problem followed by a time loop. For each time 
step, a nonlinear system is solved using the Newton--Raphson method until 
convergence (i.e. the residual $res_{\textrm{M}}$ is smaller than some
prescribed tolerance $tol_{\textrm{M}}$).
Therefore, $N_{\mathrm{GP}}$ mesoscale problems are solved in parallel and the homogenized 
law are obtained.
The term
$\Big(\frac{\bM_{\mathrm{m}}}{\Delta t_k} + \partial \bF^{(i,j,k)}_{\mathrm{m}}/\partial \bah_{\mathrm{c}}^{(i,j,k)}
\Big)^{-1} \, 
\Big(\partial \bF^{(i,j,k)}_{\mathrm{m}}/\partial \bah_{\mathrm{M}}^{(i,j,k)}\Big)$
in \eqref{eq:Jacobian-schur-complement-macro} can be interpreted as the discretization 
of the Fr\'{e}chet derivative $\Big(\partial \bB^{(i)}_{\mathrm{c}}/ \partial \bb_{\mathrm{M}} \Big)$ in 
\eqref{eq:Homogenized-Law-FE-HMM-1} (see \cite{niyonzima-chm-16}).
%
%
%
%

\subsection{The static case}
\label{sec:phd4_sec3MS_a}
The static problem can be seen as a simplified version of the dynamic problem 
obtained by neglecting the time derivatives. The macroscale weak formulation is 
derived from the $\ba-v$ formulation described in the section \ref{macro-equations}. 
The three-dimensional macroscale weak formulation reads: 
find $\ba_M \in \Hcurl[_e]{\Omega}$ such that
\begin{equation}
\left( \bh_M, \Curl[_x]{\ba_M' } \right)_{\Omega} 
= - \left< \vec{n} \times \vec{h}_M, \vec{a}_M' \right>_{\Gamma_h}
+ \left( \bj_s , \ba_M' \right)_{\Omega_s}
\label{eq:weakform_magdyn_vecPotmacro_static}
\end{equation}
holds for all test functions $\ba_M' \in \Hcurl[_e^0]{\Omega}$. 
The two-dimensional macroscale problem reads: 
find $a_{zM} \in H_e^1(\Omega)$ such that
\begin{multline}
\left( \bh_M,  \boldsymbol{1}_z \times \Grad[_x]{a_{zM}' }  \right)_{\Omega} 
\\
= - \isur[_{\Gamma_h}]{\pvec{\vec{n}}{\vec{h}_M } }{a_{zM}' \boldsymbol{1}_z} 
+ \left( j_s , a_{zM}' \right)_{\Omega_s}\,,  \qquad \forall a_{zM}' \in H_e^{10}(\Omega).
\label{eq:weakform_magdyn_vecPot_macro_2D_static}
\end{multline}

The three-dimensional mesoscale problem can be derived from
\eqref{eq:weakform_magdyn_vecPot_micro}: 
find $\ba_c \in \Hcurl[_*]{\mathcal{Y} }$ such that
\begin{equation}
  \left( \boldsymbol{\mathcal{H}}(\CurlSymb_y \ba_c + \bb_M, \bx, \by) ,  \CurlSymb_y
    \ba_c' \right)_{\Omega_{m} } = 0\,, 
  \qquad \forall \ba_c' \in \Hcurl[_*]{\mathcal{Y} }
  \label{eq:weakform_magdyn_vecPot_micro_static}
\end{equation}
and the two-dimensional mesoscale formulation reads:  
find $a_{zc} \in \Hone[_*]{\mathcal{Y} }$ such that
\begin{equation}
\left( \boldsymbol{\mathcal{H}}(\boldsymbol{1}_z \times \Grad[_y]{a_{zc}} + \bb_M, \bx, \by) , 
\boldsymbol{1}_z \times \Grad[_y]{a_{zc}' } \right)_{\Omega_{m} } = 0\,, \qquad \forall a_{zc}' \in \Hone[_*]{\mathcal{Y} }.
\label{eq:weakform_magdyn_vecPot_micro_static222}
\end{equation}

\vspace{10mm}
%
%
\section{Numerical tests}
\label{section:tests}
The models developed in the previous section are valid for the general three-dimensional 
problems. In this section we apply the models to solve nonlinear two-dimensional 
eddy current problem involving nonlinear/hysteretic materials using Problems
\ref{eq:mqs_b_conform_magdyn_a_macroscale_2D} and \ref{eq:mqs_b_conform_magdyn_a_mesoscale_2D}.

\subsection{Description of the problem}
\label{subsection:tests_descrption_b}
We consider a soft magnetic composite (SMC) material to test the ideas developed in 
the previous sections. An idealized 2D periodic SMC (with $20 \, \times \, 20$ grains) 
surrounded by an inductor is studied. We solve this academic problem using the SMC 
structure depicted in Figure \ref{fig:smc_grains_a-v} (only $10 \, \times \, 10$ grains are shown).
\begin{figure}[ht]
    \centering
    \scalebox{0.5}{\input{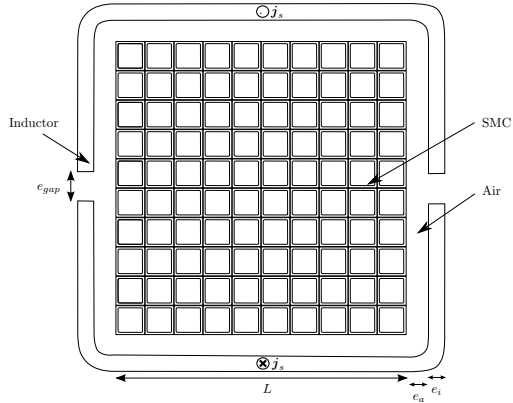}}
    \caption{Two-dimensional soft magnetic composite geometry (only 100 grains
    out of the actual 400 are drawn).
    Top and bottom inductors carry opposite source currents.
    The dimensions are $L = 1000 \, \mu$m, $e_a = 150 \, \sqrt{2}/2 \, \mu$m, 
    $e_i = 100 \, \mu$m and $e_{gap} = 100 \, \mu$m.}
    \label{fig:smc_grains_a-v}
\end{figure}
\begin{figure}[ht]
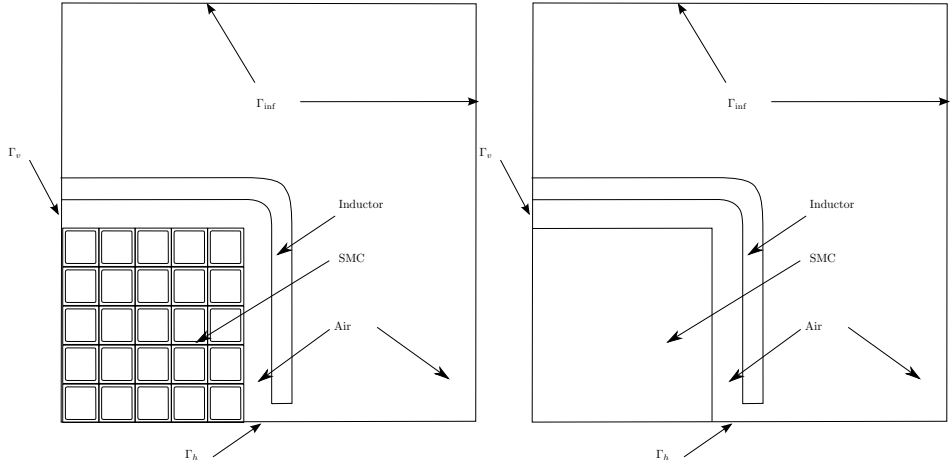

    \centering
    \hspace*{\fill}
    \scalebox{0.45}{\input{smc_multiscale_domain_2D_a-v_1_quarter.tex}}
    \hfill
    \scalebox{0.45}{\input{smc_multiscale_domain_2D_a-v_1_quarter_homog.tex}}
    \hspace*{\fill}
    \caption{Geometry used for computations (one fourth taking advantage of symmetries).  
    Left: Reference geometry (only 25 grains out of the actual 100 are depicted). 
    Right: Homogenized geometry.
    }
    \label{fig:smc_grains_a-v_quarter}
\end{figure}

The source current $\bj_s$ is imposed perpendicular to the $xy$-plane 
$\bj_s = (0, 0, j_s)$ with $j_s = j_{s0} s(t)$ where $j_{s0}$ is the amplitude 
and $s(t) = \sin(2 \pi f t)$. Therefore, the problem can be solved using 
a two-dimensional magnetic vector potential formulation with $\ba = (0, 0, a_z)$,
thus constraining the magnetic flux density $\bb$ in the $xy$-plane. 
Only one fourth the structure is considered for numerical computations thanks 
to the symmetry (see Figure~\ref{fig:smc_grains_a-v_quarter} -- Left for the 
reference geometry and Figure~\ref{fig:smc_grains_a-v_quarter} -- Right for the 
homogenized geometry).
In both cases, the following boundary conditions are imposed on 
$\Gamma_{\mathrm{inf}}, \Gamma_{h}$ and $\Gamma_{v}$:
\begin{align}
(\bn \cdot \bb) |_{\Gamma_{\mathrm{inf}}} = 0 
&\quad \Leftarrow \quad 
(\bn \times \ba) |_{\Gamma_{\mathrm{inf}}} = \boldsymbol{0},
\\
(\bn \cdot \bb) |_{\Gamma_{\mathrm{h}}} = 0 
&\quad \Leftarrow \quad 
(\bn \times \ba) |_{\Gamma_{\mathrm{h}}} = \boldsymbol{0}, \quad (\bn \times \bh) |_{\Gamma_{\mathrm{v}}} = \boldsymbol{0}.
\label{eq:chap_results_BC}
\end{align}

We consider operating frequencies smaller than 50\,kHz, which 
corresponds to $\lambda_f \text{ and } \lambda_M \simeq \lambda = 6000$m). The smallest 
wavelength of the source is much larger than the length of the structure ($\simeq 500 \mu m$) 
so that the  assumption of a magnetoquasistatic problem can be made.

All materials are isotropic, so that the magnetic field $\bh$ has only 
$xy$ components. The conducting grains (electric conductivity $\sigma = 5 \,
10^{6}$\,S/m) are surrounded by a perfect insulator, linear and
non-magnetic ($\mu_r = 1$). The grains are governed by the following magnetic laws:
\begin{enumerate}
\item a nonlinear exponential law $\boldsymbol{\mathcal{H}}(\bb) 
= \left( \alpha + \beta \, \exp (\gamma ||\bb||^2) \right) \, \bb$ with 
$\alpha=  388, \beta = 0.3774$ and $\gamma = 2.97$ \cite{delince-thesis-94}.
\item a Jiles--Atherton hysteresis model with parameters 
$\boldsymbol{\mathcal{M}}_s = 1,145,500$\,A/m, $a = 59$\,A/m, $k = 99$\,A/m, $c = 0.55$ and $\alpha = 1.3 \, 10^{-4}$ 
(see \cite{gyselinck-hyst-04,benabou-hyst-03} for more details on the Jiles--Atherton model and 
the meaning of the parameters it uses).  
\end{enumerate}

Results obtained using the computational homogenization (subscript ``comp'' 
for \emph{computational homogenization}, subscript ``M'' for \emph{Macro} and ``m'' for \emph{meso})
are compared to the reference results (subscript ``Ref'') obtained solving the
reference problem (i.e. the weak form of
(\ref{eq:equations_PE_multiscale} a-b)--(\ref{eq:materiallaw_PE_multiscale} a-b) on
a very fine mesh.

Some quantities of interest (global quantities and errors) are defined and used 
for numerical validation. The global quantities are the reference and the 
computational homogenization eddy currents losses:
\begin{equation}
\begin{aligned}
&\tau \textrm{P}_{\mathrm{Ref}}(t) = \displaystyle \int_{\Omega_c} (\sigma |\partial_t \ba^{\varepsilon}(\bx, t)|^2) \, \mathrm{d} x, 
\\
&\tau \textrm{P}_{\mathrm{m}}(t) = \displaystyle \int_{\Omega} \tau \textrm{P}_{\mathrm{m}}^{\mathrm{up}}(\bx, t) \, \mathrm{d} x
= 
\displaystyle \int_{\Omega} \Big( \frac{1}{|\Omega_m|} \int_{\Omega_{\mathrm{mc}}} 
(\sigma |\partial_t \ba_{\mathrm{m}}(\bx, \by, t)|^2) \, \mathrm{d} y \Big) \, \mathrm{d} x.
\label{eq:ecl_all}
\end{aligned}
\end{equation}
The equivalent quantities in terms of the magnetic energy can be defined.

Two types of errors are defined as:
\begin{itemize}
\item the relative error in terms of the eddy current losses:  
\begin{equation}
\mathrm{Err}_{\tau \textrm{P}} 
= \dfrac{\parallel \tau \textrm{P}_{\mathrm{m}}
- \tau \textrm{P}_{\mathrm{Ref}}\parallel_{L^{\infty}(0,\, T)}  }
{\parallel\tau \textrm{P}_{\mathrm{Ref}}\parallel_{L^{\infty}(0,\, T)}},
\label{eq:error-ECL}
\end{equation}
\item the pointwise relative error on the fields $\bb_\mathrm{M}$ and $\bb_\mathrm{m}$:
\begin{multline}
\mathrm{Err}_{\mathrm{M}}^{\bb}(\bx) =
\dfrac{\parallel \bb_{\mathrm{M}}(\bx) - \bb_{\mathrm{Ref}}(\bx) \parallel_{\Ltwo{0,\, T} } }
{\parallel\bb_{\mathrm{Ref}}(\bx) \parallel_{\Ltwo{0,\, T} } },
\,
\\
\mathrm{Err}_{\mathrm{m}}^{\bb}(\bx) =  
\dfrac{\parallel \bb_{\mathrm{m}}(\bx) - \bb_{\mathrm{Ref}}(\bx) \parallel_{\Ltwo{0,\, T} } }
{\parallel\bb_{\mathrm{Ref}}(\bx) \parallel_{\Ltwo{0,\, T} } },
\label{eq:error-fields-b}
\end{multline}
\end{itemize}

\subsection{Results} 
\label{sec:phd_4_smc_b_results}
Results of the reference and the multiscale problems are compared in this section.
The latter are obtained by solving a finite element problem on the entire, 
finely meshed multiscale domain (110,282 triangular elements). 
%
%
Computational results are carried out on a macroscale, coarse mesh (42 quad elements). 
Mesoscale problems are solved around each numerical quadrature point of the
macroscale mesh using a fine mesh (4125 triangular elements). 

\begin{figure}[h]
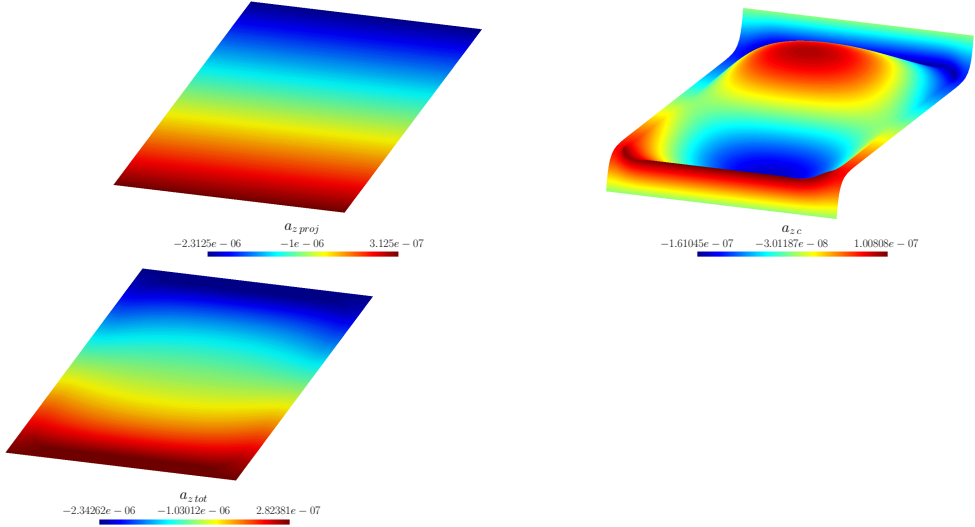

\centering
\hspace*{\fill}
\scalebox{0.12}{\input{a_proj_p1_cs3_gp899994.tex}}
\hfill
\scalebox{0.12}{\input{a_pert_p1_cs3_gp899994.tex}}
\hfill
\scalebox{0.12}{\input{a_tot_p1_cs3_gp899994.tex}}
\hspace*{\fill}
\caption{Terms contributing to the total mesoscale magnetic vector potential
  for a cell problem centered at ($325, 25, 0.0$)$\mu m$. 
Top: the $z$-component of the projection term $\ba_{proj}(\bx, \by, t) =
\ba_M(\bx, t) + \kappa (\by \times \bb_M(\bx, t))$. 
Middle: the $z$-component of the correction term $\ba_c(\bx, \by, t)$. 
Bottom: the $z$-component of the total mesoscale vector potential $\ba_{tot}(\bx, \by, t)$ \{nonlinear case with $j_{s0} = 35 \times 10^{7} \text{A/m}^2, f = 25\text{ kHz}$\}. }
\label{fig:smc_a-v_mesofields}
\end{figure}
Figure~\ref{fig:smc_a-v_mesofields} depicts the different contributing terms 
involved in the resolution of the mesoscale problem. The projection term which 
varies linearly on the mesoscale domain is computed from the macroscale fields 
as $\ba_{proj}(\bx, \by, t) = \ba_M(\bx, t) + \kappa (\by \times \bb_M(\bx, t))$. 
This term is then used as a source for the computation of the correction term 
$\ba_{c}(\bx, \by, t)$ at the mesoscale level which allows to derive the total 
magnetic vector potential $\ba_{tot}(\bx, \by, t) = \ba_c(\bx, \by, t) + \ba_M(\bx, t) + \kappa (\by \times \bb_M(\bx, t))$.
%
%
%
%
\begin{figure}[h]
\centering

  \includegraphics[width=0.5\textwidth]{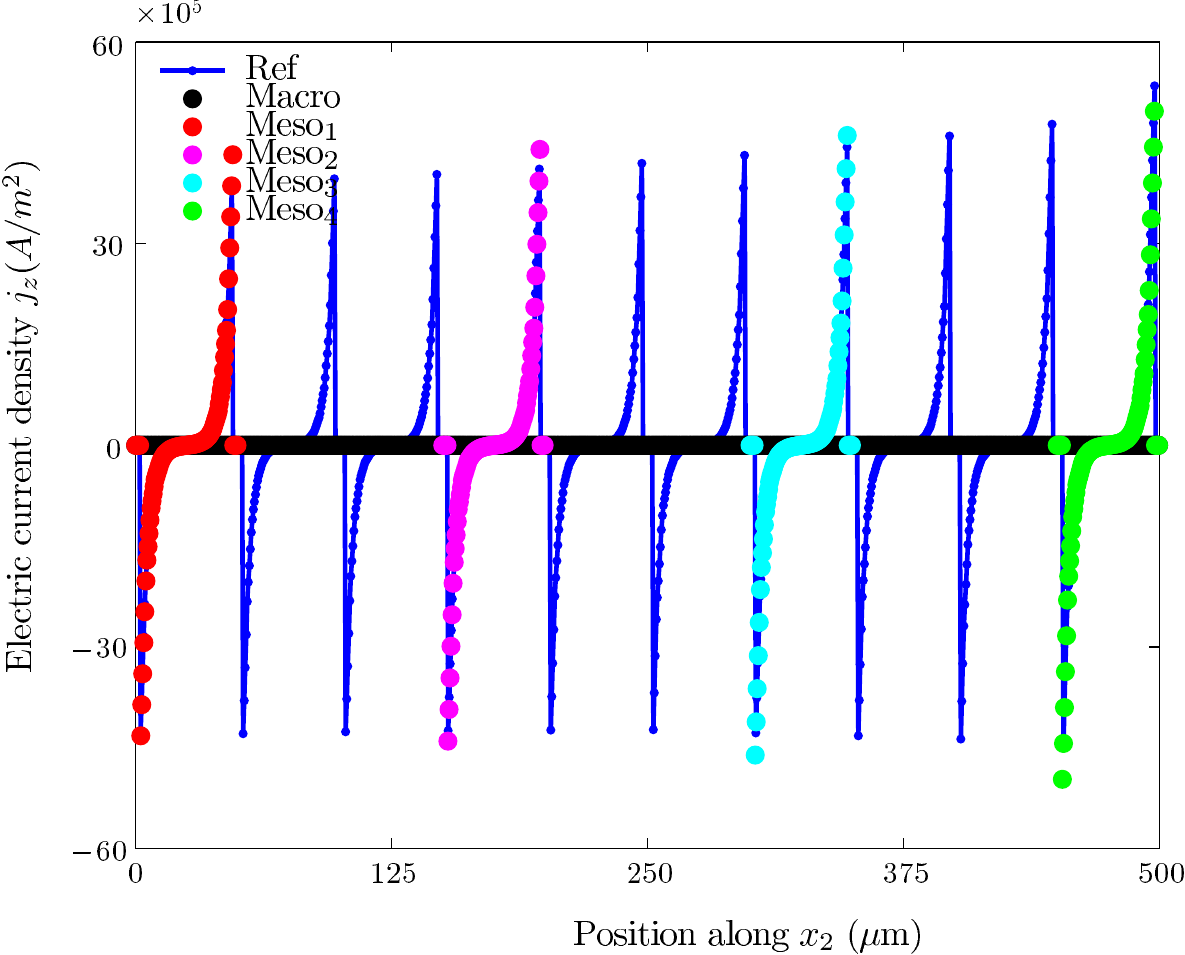}
  \includegraphics[width=0.5\textwidth]{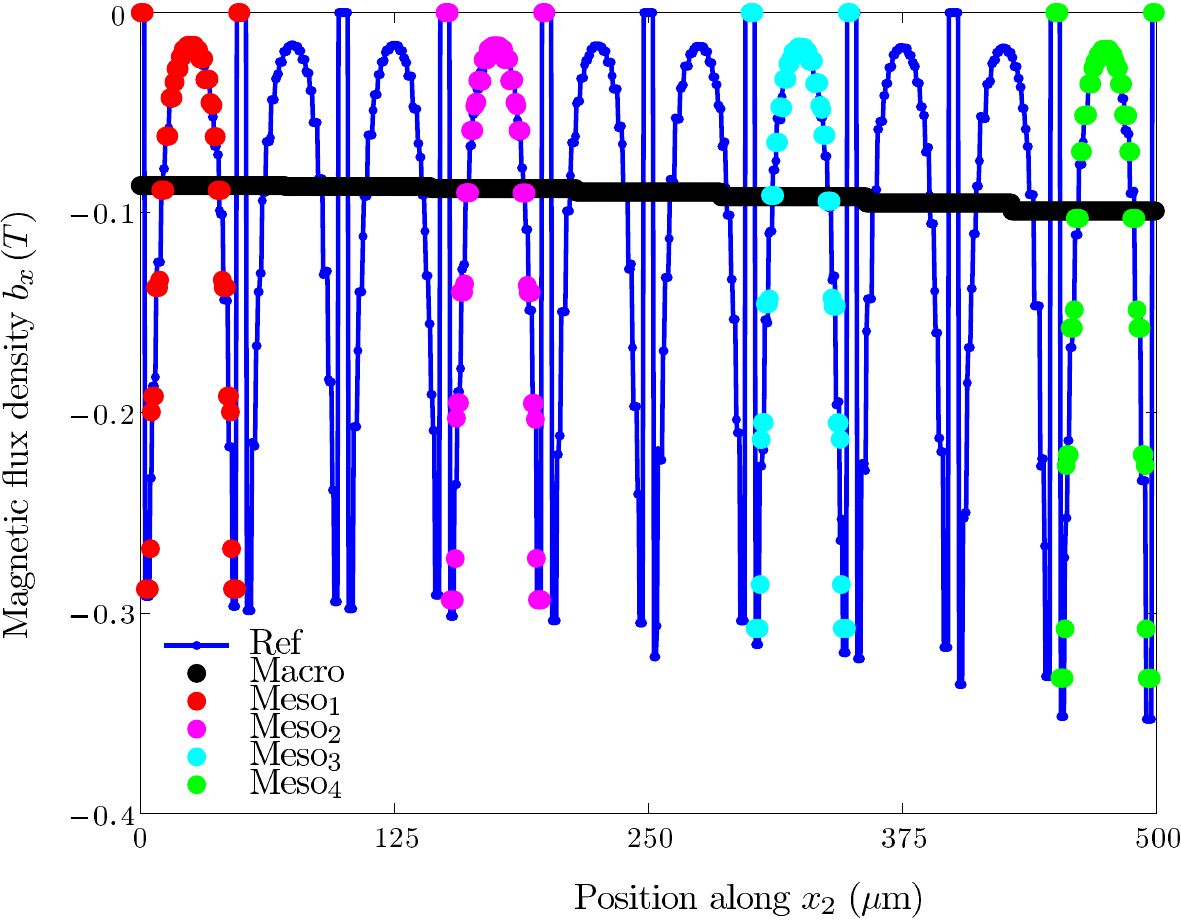}
  \includegraphics[width=0.5\textwidth]{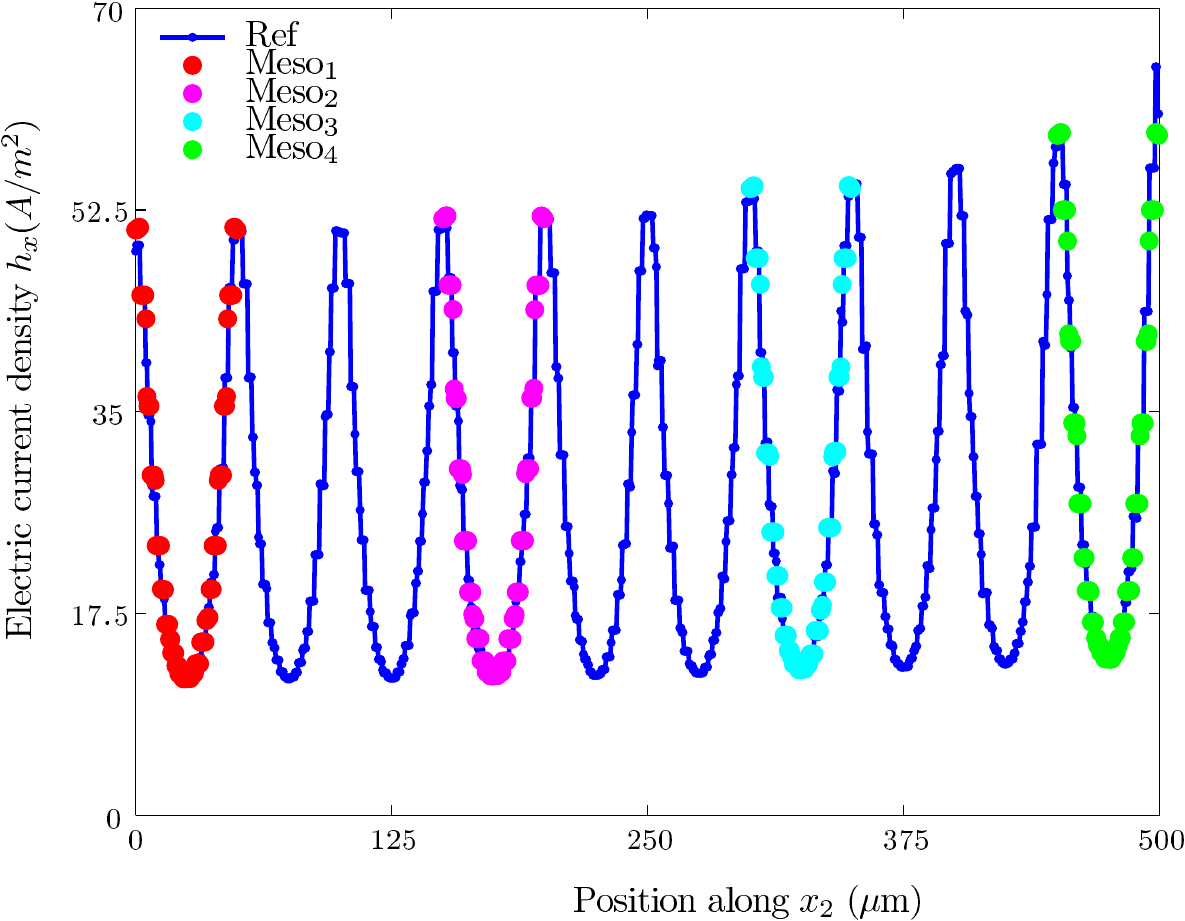}
\caption{SMC problem, $\vec{b}$-conform formulations, hysteretic case. 
Spatial cuts of the $z$-component of the eddy currents $\vec{j}$ (top), of 
the $x$-component of the magnetic induction $\vec{b}$ (middle) and of the magnetic field $\vec{h}$ (bottom) along the line 
\{$x = 25, z = 0$\}\,$\mu$m. ($f = 10$\,kHz, $t = 5 \, 10^{-7}$s for the curve 
of eddy currents and $t = 25 \, 10^{-7}$s for the curve of the magnetic induction).} 
\label{fig:smc_a-v_local_cuts_hysteresis}
\end{figure}

The spatial cuts of the magnetic induction $\bb$, the eddy 
currents $\bj$ and the magnetic field $bh$ are shown in~Figure
\ref{fig:smc_a-v_local_cuts_hysteresis}. The agreement between the
reference solution and the mesoscale solution on a cell around certain
Gau\ss \, points in the computational domain proves excellent. 
As expected, small discrepancies are observed near the boundary of the domain 
(see Tables~\ref{tab:b-conform_localcuts_b_hyst} and~\ref{tab:b-conform_L2_error_b_hyst}).

Table~\ref{tab:b-conform_localcuts_b_hyst} displays the values of $\parallel\bb\parallel$ 
obtained from the reference solution (Reference), the macroscale solution (Macro) 
and the mesoscale solution (Meso) and the corresponding relative pointwise errors 
(Error meso, Error macro) for $t = 6 \cdot 10^{-6}$\,s.
In this table, we observe that the mesoscale error increases with the
proximity to the boundary of the computational domain. In the bulk, the error is around 1\% and rises up
to 14\% at the boundary.
Indeed, cells located near the boundary do not respect the periodicity
assumption, they are not immersed in a periodic environment.
The macroscale error is huge and almost independent of the location of the considered point.
\begin{table}[t]
\caption{
Soft magnetic composite problem -- $\vec{b}$-conform formulations. 
Comparison of the reference and the computational 
(macroscale and mesoscale) magnetic flux density ($\|\vec{b}\|$\,[T]) 
at different points of the macroscale domain \{$t = 6 \cdot 10^{-6}$\,s\}.
}
\label{tab:b-conform_localcuts_b_hyst}
\medskip\centering
\begin{tabular}{| c | c | c | c | c | c |}
\hline
Position ($\mu$m)   &  Reference   & Meso       &  Macro    & $\mathrm{Err}_{\mathrm{m}}^{\bb}(\bx) (\%)$ & $\mathrm{Err}_{\mathrm{M}}^{\bb}(\bx) (\%)$\\
\hline
 $(25, 25, 0)$      &  0.0157652   & 0.0158937  & 0.0347775 & 0.82       & 120.60\\
 $(25, 475, 0)$     &  0.0186482   & 0.0181317  & 0.0403767 & 2.77       & 116.52\\
 $(175, 175, 0)$    &  0.0158077   & 0.0158738  & 0.0346577 & 0.42       & 119.25\\
 $(475, 25, 0)$     &  0.0156693   & 0.0158615  & 0.0345838 & 1.23       & 120.70\\
 $(475, 475, 0)$    &  0.0184396   & 0.0158563  & 0.0417285 & 14.01      & 126.30\\
\hline
\end{tabular}
\end{table}

Table~\ref{tab:b-conform_L2_error_b_hyst} provides the relative $L^2(0, T)$
error defined by~\eqref{eq:error-fields-b}.
Results of Table~\ref{tab:b-conform_L2_error_b_hyst} allow us to draw the
same conclusions as those from Table~\ref{tab:b-conform_localcuts_b_hyst}, i.e.\ the error increases 
as the point gets closer to the boundary of the computational domain. 
\begin{table}[t]
\caption{Soft magnetic composite problem -- $\vec{b}$-conform formulations. 
Relative $L^{2}(0, T)$ errors of the mesoscale and the macroscale magnetic
flux density with regard to the reference,
$\mathrm{Err}_{\mathrm{m}}^{\bb}(\bx)$ and
$\mathrm{Err}_{\mathrm{M}}^{\bb}(\bx)$, respectively, at different points of
the computational domain.
}
\label{tab:b-conform_L2_error_b_hyst}
\medskip\centering
\begin{tabular}{| c | c | c |}
\hline
Position ($\mu m$)  & $\mathrm{Err}_{\mathrm{m}}^{\bb}(\bx) (\%)$ & $\mathrm{Err}_{\mathrm{M}}^{\bb}(\bx) (\%)$ \\
\hline
 $(25, 25, 0)$      & 3.27                                   & 11.49                                  \\
 $(25, 475, 0)$     & 4.93                                   & 15.13                                  \\
 $(175, 175, 0)$    & 3.01                                   & 11.88                                  \\
 $(475, 25, 0)$     & 3.04                                   & 12.27                                  \\
 $(475, 475, 0)$    & 15.46                                  & 22.91                                  \\
\hline
\end{tabular}
\end{table}
%
%

%
%
\begin{figure}[ht]
\centering
\Huge
\resizebox{0.95\textwidth}{!}{\input{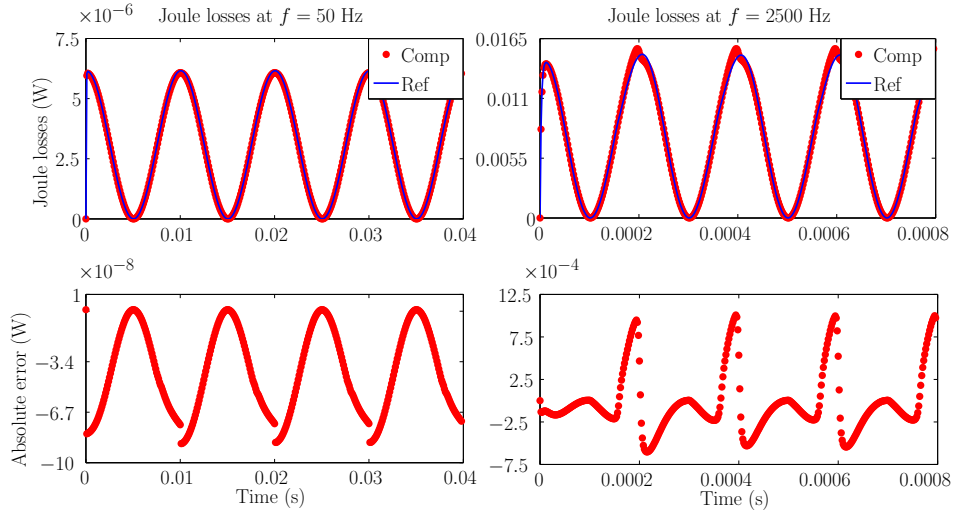}}
\caption{SMC problem, $\vec{b}$-conform formulations, hysteretic case. 
Instantaneous Joule losses and absolute error between the reference (Ref) and 
the computational  (Comp) solutions. Hysteretic case. Left: $f = 50$\,Hz. Right: $f = 2500$\,Hz.}
\label{fig:smc_a-v_Joule-losses}
\end{figure}
Figure \ref{fig:smc_a-v_Joule-losses} depicts the evolution of the eddy currents losses
for excitations at $50$\,Hz and $2500$\,Hz (which correspond to 
the case with enhanced skin effect). A good agreement between Joules losses is 
observed for both frequencies: a maximum error of $1.41\,\%$ and $6.69\,\%$ 
are observed for $f = 50$\,Hz and $f = 2500$\,Hz, respectively. 

Table~\ref{tab:b-conform_error_freq_Joulelosses} contains the relative $L^{\infty}(0, T)$ 
error of the Joule losses defined by equation \eqref{eq:error-ECL} as a function of frequency. 
\begin{table}[ht]
\caption{Soft magnetic composite problem -- $\vec{b}$-conform formulations. 
Relative $L^{\infty}(0,T)$ norm error on the total Joule losses as a function of the frequency.}
\label{tab:b-conform_error_freq_Joulelosses}
\medskip\centering
\begin{tabular}{| c | c | c |}
\hline
Frequency (Hz)         & $\mathrm{Err}_{\tau \textrm{P}}$  ($\%$) \\
\hline
 $50$                  & $1.41$                                 \\
 $100$                 & $1.46$                                 \\
 $250$                 & $1.61$                                 \\
 $1000$               & $3.42$                                 \\
 $2500$               & $6.69$                                 \\
\hline
\end{tabular}
\end{table}

Figure \ref{fig:smc_a-v_hyst_convergence} shows the convergence of the 
residual resulting from the resolution by the Newton--Raphson method as a function 
of the number of nonlinear iteration. It can be seen that the macroscale problem 
converges quadratically while the mesoscale problems converge at an average rate of 1.33.
\begin{figure}[ht]
\begin{center} 
\hspace*{\fill} 
\includegraphics[width=0.47\textwidth]{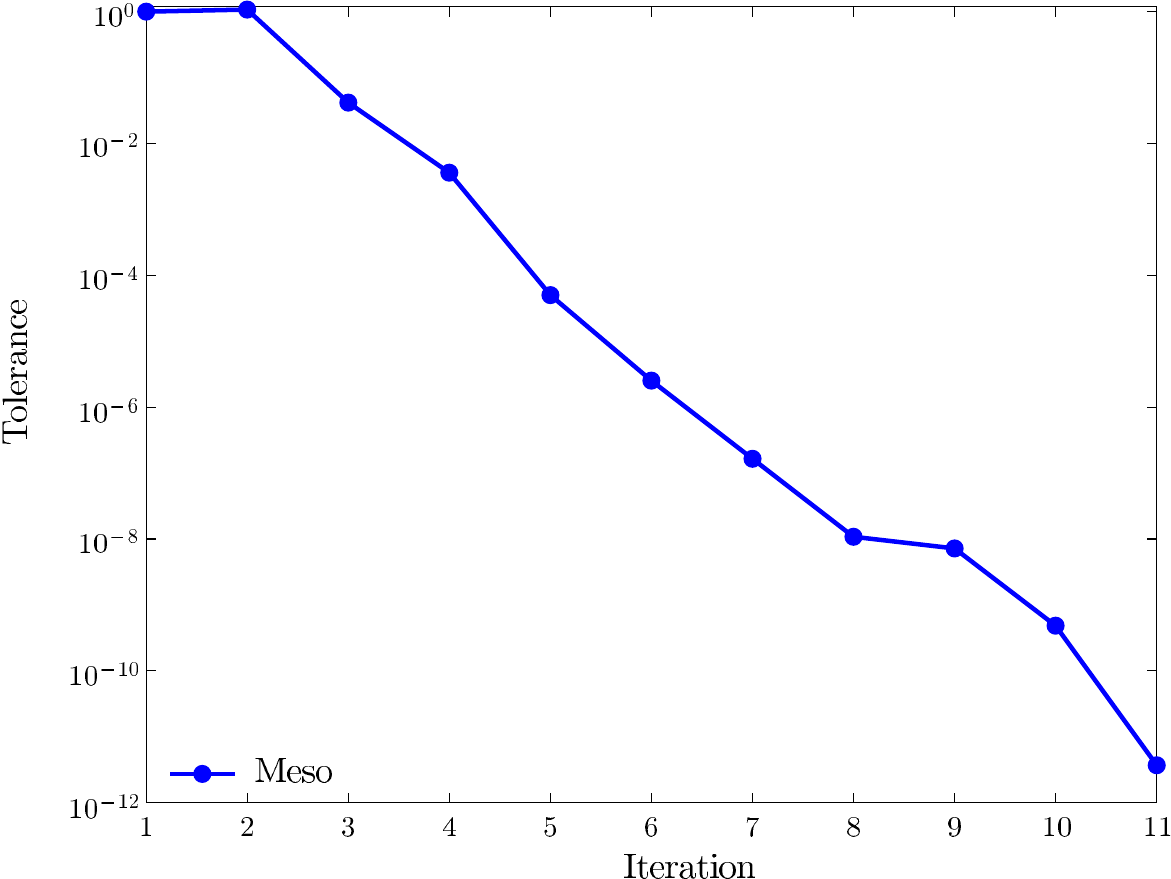}
\hfill 
\includegraphics[width=0.47\textwidth]{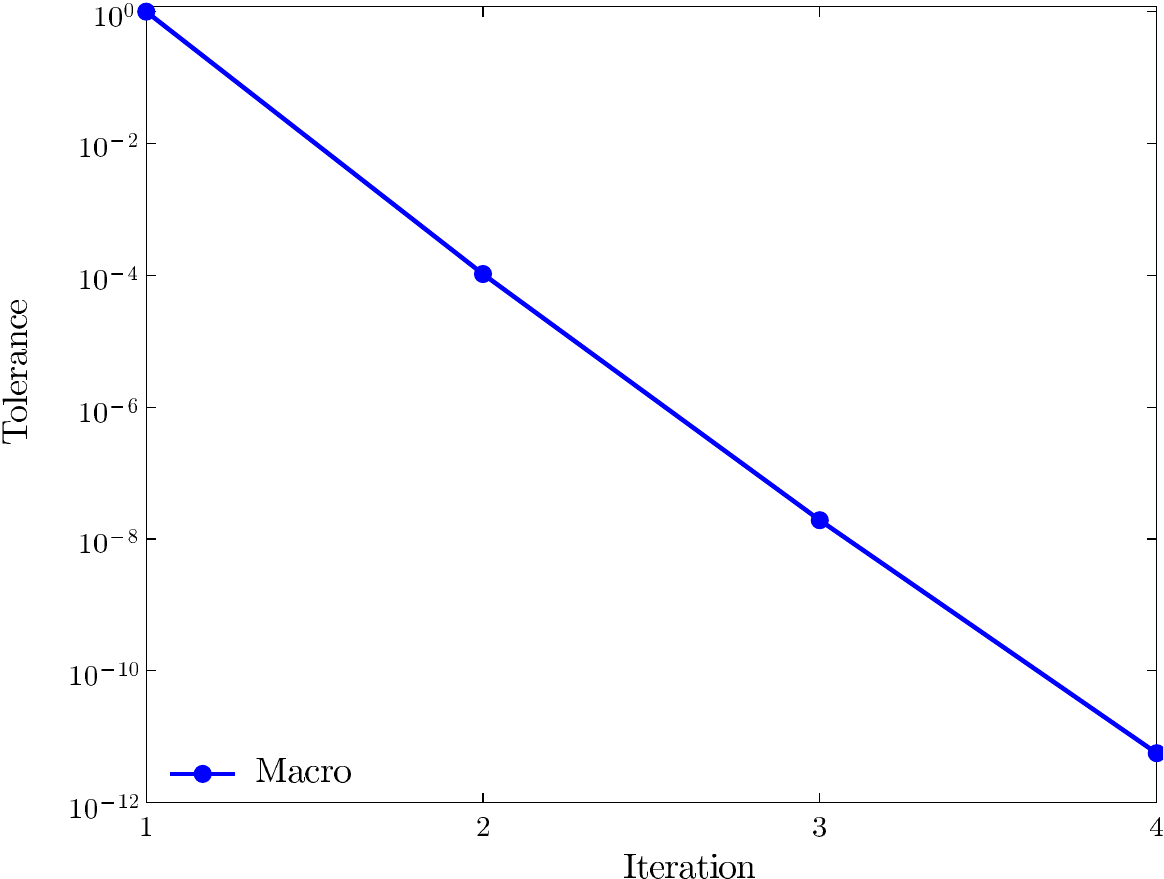}
\hspace*{\fill} 
\end{center}
\caption{SMC problem, $\vec{b}$-conform formulations, hysteretic case. Convergence of the error as a function of nonlinear iterations. Top: mesoscale problem. Bottom: macroscale problem.}
\label{fig:smc_a-v_hyst_convergence}
\end{figure}

\section{Conclusions}
\label{section:conclusions}

In this paper we have developed a computational multiscale method inspired by the 
HMM approach to solve nonlinear, possibly hysteretic magnetoquasistatic problems on multiscale domains (e.g. composite 
materials, lamination stacks, etc.). 
To construct the computational multiscale model, we combine theoretical 
results from two-scale convergence theory and asymptotic homogenization. 
The two-scale convergence and periodic unfolding methods are used for deriving the partial differential 
equations governing fields at both the macroscale and the mesoscale levels, 
valid in the nonlinear regime and in the presence of $\CurlSymb$ differential operators. 
Asymptotic homogenization is used for defining a mesoscale problem in the case of linear 
constitutive laws (e.g. the linear electric conductivity law). 

Although the theoretical foundation is only valid in the case of linear and 
nonlinear problems governed by a maximal monotone operator, in practice, the 
resulting numerical multiscale scheme has been successfully applied to general 
magnetoquasistatic problems also exhibiting memory effects (hysteresis). 
The numerical tests were performed for magnetodynamic problems, 
using $\bb$-conform formulations. An excellent agreement has been obtained between 
the reference solutions (computed using a brute force approach) and the computational 
(mesoscale) solutions. 
\textcolor{black}{We observed larger errors} near the boundary
of the computational domain as the cell problems defined near the boundary are not immersed 
in a periodic environment. The eddy current losses are also accurately evaluated. 
The error on these losses increases as a function of the frequency.

\textcolor{black}{For the considered academic test case, the proposed computational multiscale method 
fulfills the original goals (Section \ref{sec:introduction}): it allows to solve 
multiscale magnetoquasistatic problems, including the computation of local fields at 
the mesoscale and the accurate evaluation of electromagnetic losses.
It naturally handles nonlinear or hysteretic materials and periodic mesoscale 
geometries. From an engineering point of view, the approach could be straightforwadly applied 
to deal with more complex multiscale geometries.}

The main disadvantage of the method is its higher computational cost. 
However, since all the mesoscale problems are independent, the method is perfectly suited 
for modern massively parallel computers, and we thus believe that it has a lot 
of potential, even compared to brute force approaches, which do not scale well.

\section*{Acknowledgment}
This work was supported by the the Belgian Science Policy under grant IAP P7/02 
(Multiscale modelling of electrical energy system). Patrick Dular is a fellow 
with the Fonds de la recherche scientifique-FNRS (FRS-FNRS).


\begin{thebibliography}{10}

\bibitem{Abdulle:2009p1255}
{\sc A.~Abdulle}, {\em The finite element heterogeneous multiscale method: a
  computational strategy for multiscale {PDEs}}, GAKUTO International Series
  Math. Sci. Appl., Multiple scales problems in Biomathematics, Mechanics,
  Physics and Numerics, 31 (2009), pp.~133--181.

\bibitem{Abdulle:2003p808}
{\sc A.~Abdulle and W.~E}, {\em Finite difference heterogeneous multi-scale
  method for homogenization problems}, Journal of Computational Physics, 191
  (2003), pp.~18--39,
  \href{http://dx.doi.org/10.1016/S0021-9991(03)00303-6}{doi:\nolinkurl{10.1016/S0021-9991(03)00303-6}}.

\bibitem{acevedo-mqs-13}
{\sc R.~Acevedo and G.~Loaiza}, {\em A fully-discrete finite element
  approximation for the eddy currents problem}, Ingenier{\'\i}a y Ciencia, 9
  (2013), pp.~111--145.

\bibitem{bachinger-mqs-05}
{\sc F.~Bachinger, U.~Langer, and J.~Sch{\"o}berl}, {\em Numerical analysis of
  nonlinear multiharmonic eddy current problems}, Numerische Mathematik, 100
  (2005), pp.~593--616.

\bibitem{belkadi-homogenization-09}
{\sc M.~Belkadi, B.~Ramdane, D.~Trichet, and J.~Fouladgar}, {\em Non linear
  homogenization for calculation of electromagnetic properties of soft magnetic
  composite materials}, IEEE: Transaction on Magnetics, 45 (2009),
  pp.~4317--4320.

\bibitem{benabou-hyst-03}
{\sc A.~Benabou, S.~Cl\'{e}net, and F.~Piriou}, {\em Comparison of {P}reisach
  and {J}iles-{A}therton models to take into account hysteresis phenomenon for
  finite element analysis}, Journal of Magnetism and Magnetic Materials, 261
  (2003), pp.~305--310.

\bibitem{bensoussan-ahm-11}
{\sc A.~Bensoussan, J.-L. Lions, and G.~Papanicolaou}, {\em Asymptotic Analysis
  for Periodic Structures}, American Mathematical Society, 2011.

\bibitem{bossavit-modelisation-93}
{\sc A.~Bossavit}, {\em {\'E}lectromagn\'etisme, en vue de la mod\'elisation},
  Springer-Verlag, 1993.

\bibitem{Bossavit1994}
{\sc A.~Bossavit}, {\em Effective penetration depth in spatially periodic
  grids: a novel approach to homogenization}, in Proceedings, 1994,
  pp.~859--864.

\bibitem{Bossavit1996}
{\sc A.~Bossavit}, {\em Homogenizing spatially periodic materials with respect
  to maxwell equations: Chiral materials by mixing simple ones}, in
  Proceedings, 1996, pp.~564--567.

\bibitem{bossavit-computational-98}
{\sc A.~Bossavit}, {\em Computational Electromagnetism. {V}ariational
  Formulations, Edge Elements, Complementarity}, Academic Press, 1998.

\bibitem{bottauscio-chm-08}
{\sc O.~Bottauscio, V.~Chiado~Piat, M.~Chiampi, M.~Codegone, and A.~Manzin},
  {\em Nonlinear homogenization technique for saturable soft magnetic
  composites}, IEEE Transactions on Magnetics, 44 (2008), pp.~2955--2958.

\bibitem{bottauscio-msfem-13-1}
{\sc O.~Bottauscio, M.~Chiampi, and A.~Manzin}, {\em Multiscale modeling of
  heterogeneous magnetic materials}, International Journal of numerical
  modeling: electronic networks, devices and fields, 27 (2014), pp.~373--384.

\bibitem{bottauscio-msfem-13}
{\sc O.~Bottauscio and A.~Manzin}, {\em Comparison of multiscale models for
  eddy current computation in granular magnetic materials}, Journal of
  Computational Physics, 253 (2013), pp.~1--17.

\bibitem{braides-gamma-02}
{\sc A.~Braides}, {\em {$\Gamma$}-convergence for Beginners}, vol.~22,
  Clarendon Press, 2002.

\bibitem{laurence-mfh-10}
{\sc L.~Brassart, I.~Doghri, and D.~L.}, {\em Homogenization of elasto-plastic
  composites coupled with a nonlinear finite element analysis of the equivalent
  inclusion problem}, International Journal of Solids and Structures, 47
  (2010), pp.~716--729.

\bibitem{brezzi-vms-97}
{\sc F.~Brezzi, L.~P. Franca, T.~J.~R. Hughes, and A.~Russo}, {\em $b = \int
  g$}, Computer Methods in Applied Mechanics and Engineering, 145 (1997),
  pp.~329--339.

\bibitem{cioranescu-pum-02}
{\sc D.~Cioranescu, A.~Damlamian, and G.~Griso}, {\em Periodic unfolding and
  homogenization}, C.R. Acad. Sci. Paris, Ser. I, 335 (2002), pp.~99--104.

\bibitem{cioranescu-pum-08}
{\sc D.~Cioranescu, P.~Donato, and R.~Zaki}, {\em The periodic unfolding method
  in homogenization}, S.I.A.M. J. Math. Anal., 40 (2008), pp.~1585--1620.

\bibitem{corcolle-thesis-09}
{\sc R.~Corcolle}, {\em D\'{e}termination de Lois de Comportment Coupl\'{e} par
  des Techniques d'Homog\'{e}n\'{e}isation: application aux Mat\'{e}riaux du
  G\'{e}nie Electrique}, PhD thesis, Universite Paris-Sud XI, 2009.

\bibitem{DalMaso-gammaconv-93}
{\sc G.~Dal~Maso}, {\em Introduction to {$\Gamma$}-Convergence}, Birkhauser,
  1993.

\bibitem{DeGiorgio-gammaconv-84}
{\sc E.~De~Giorgi}, {\em G-operators and {$\Gamma$}-convergence}, In Proc. Int.
  Congr. Math.,  (1984), pp.~1175--1191.

\bibitem{delince-thesis-94}
{\sc F.~Delinc\'{e}}, {\em Mod\'elisation des R\'egimes Transitoires dans les
  Syst\`{e}mes Comportant des Mat\'{e}riaux Magn\'{e}tiques Non-Lin\'{e}aires
  et Hyst\'{e}r\'{e}tiques}, PhD thesis, Universit\'{e} de Li\`{e}ge, 1994.

\bibitem{deriaz-helmholtzdecomposition-09}
{\sc E.~Deriaz and V.~Perrier}, {\em {Orthogonal Helmholtz decomposition in
  arbitrary dimension using divergence-free and curl-free wavelets}}, Applied
  and Computational Harmonic Analysis, 26 (2009), pp.~249--269.

\bibitem{00.compumag99.dular.MagDynHA}
{\sc P.~Dular, P.~Kuo-Peng, C.~Geuzaine, N.~Sadowski, and J.~P.~A. Bastos},
  {\em Dual magnetodynamic formulations and their source fields associated with
  massive and stranded inductors}, IEEE Transactions on Magnetics, 36 (2000),
  pp.~1293--1299.

\bibitem{E:2003p1296}
{\sc W.~E}, {\em Analysis of the heterogeneous multiscale method for ordinary
  differential equations}, Comm. Math. Sci., 1 (2003), pp.~423--436.

\bibitem{e-hmm-11}
{\sc W.~E}, {\em Principles of Multiscale Modeling}, Cambridge, 2011.

\bibitem{E:2003p1295}
{\sc W.~E and B.~Engquist}, {\em The heterogeneous multiscale methods}, Comm.
  Math. Sci., 1 (2003), pp.~87--132.

\bibitem{E:2003p1294}
{\sc W.~E and B.~Engquist}, {\em Multiscale modeling and computation}, Notices
  Amer. Math. Soc., 50 (2003), pp.~1062--1070.

\bibitem{E:2003p1288}
{\sc W.~E, B.~Engquist, and Z.~Huang}, {\em Heterogeneous multiscale method: A
  general methodology for multiscale modeling}, Physical Review B, 67 (2003),
  p.~092101,
  \href{http://dx.doi.org/10.1103/PhysRevB.67.092101}{doi:\nolinkurl{10.1103/PhysRevB.67.092101}}.

\bibitem{e-hmm-07}
{\sc W.~E, B.~Engquist, X.~Li, W.~Ren, and E.~Vanden-Eijnden}, {\em
  Heterogeneous multiscale methods: A review}, Communications in Computational
  Physics, 3 (2007), pp.~367--450.

\bibitem{efendiev-msfem-04-1}
{\sc Y.~Efendiev, T.~Hou, and V.~Ginting}, {\em Multiscale finite element
  methods for nonlinear partial differential equations}, Communications in
  Mathematical Sciences, 2 (2004), pp.~553--589.

\bibitem{ekeland-convexanalysis-69}
{\sc I.~Ekeland and R.~Temam}, {\em Analyse convexe et probl\`{e}mes
  variationnelles}, Dunod, Paris, 1974.

\bibitem{feddi-thm-97}
{\sc M.~El~Feddi, Z.~Ren, A.~Razek, and A.~Bossavit}, {\em Homogenization
  technique for {M}axwell equations in periodic structure}, IEEE Transactions
  on Magnetics, 33 (1997), pp.~1382--1385.

\bibitem{evans-pde-10}
{\sc L.~C. Evans}, {\em Partial Differential Equations}, American Mathematical
  Society, Providence, Rhode Island, 2010.

\bibitem{geers-fe2-01}
{\sc M.~G.~D. Geers, V.~G. Kouznetsova, and Brekelmans}, {\em Gradient-enhanced
  computational homogenization for the micro-macro scale transition}, Journal
  de Physique IV, 11 (2001), pp.~5145--5152.

\bibitem{gyselinck-hyst-04}
{\sc J.~Gyselinck}, {\em Incorporation of a {J}iles-{A}therton vector
  hysteresis model in {2-D} {FE} magnetic computations}, COMPEL: The
  International Journal for Computation and Mathematics in Electrical and
  Electronic Engineering, 23 (2004), pp.~685--693.

\bibitem{gyselinck-homogenization-04}
{\sc J.~Gyselinck and P.~Dular}, {\em A time-domain homogenization technique
  for laminated iron cores in {3-D} finite element models}, IEEE: Transaction
  on Magnetics, 40 (2004), pp.~856--859.

\bibitem{gyselinck-homogenization-06}
{\sc J.~Gyselinck, R.~V. Sabariego, and P.~Dular}, {\em A nonlinear time-domain
  homogenization technique for laminated iron cores in three-dimensional finite
  element models}, IEEE: Transaction on Magnetics, 42 (2006), pp.~763--766.

\bibitem{harouna-helmholtzdecomposition-10}
{\sc S.~K. Harouna and V.~Perrier}, {\em {Helmholtz-Hodge Decomposition on [0,
  1]$^d$ by Divergence-free and Curl-free Wavelets}}, in International
  Conference on Curves and Surfaces, Springer, 2010, pp.~311--329.

\bibitem{hou-msfem-97}
{\sc T.~Y. Hou and X.~H. Wu}, {\em A multiscale finite element method for
  elliptic problems in composite materials and porous media}, Journal of
  Computational Physics, 134 (1997), pp.~169--189.

\bibitem{jackson1999classical}
{\sc J.~D. Jackson}, {\em Classical electrodynamics}, Wiley, 1999.

\bibitem{jiang-mqs-12}
{\sc X.~Jiang and W.~Zheng}, {\em An efficient eddy current model for nonlinear
  maxwell equations with laminated conductors}, SIAM Journal on Applied
  Mathematics, 72 (2012), pp.~1021--1040.

\bibitem{juanes-vms-05}
{\sc R.~Juanes and T.~W. Patzek}, {\em A variational multiscale finite element
  method for multiphase flow in porous media}, Finite Elements in Analysis and
  Design, 41 (2005), pp.~763--777.

\bibitem{kouznetsova-fe2-01}
{\sc V.~G. Kouznetsova, W.~A.~M. Brekelmans, and F.~P.~T. Baaijens}, {\em An
  approach to micro-macro modeling of heterogeneous materials}, Computational
  Mechanics, 27 (2001), pp.~37--48.

\bibitem{ledger-formulations-10}
{\sc P.~Ledger and S.~Zaglmayr}, {\em hp-finite element simulation of
  three-dimensional eddy current problems on multiply connected domains},
  Computer Methods in Applied Mechanics and Engineering, 199 (2010),
  pp.~3386--3401.

\bibitem{maxwell-garnett-mixing-04}
{\sc J.~C. Maxwell~Garnett}, {\em Colors in metal glasses and metal films},
  Phil. Trans. R. Soc. Lond. A, 203 (1904), pp.~385--420.

\bibitem{meunier-chm-10}
{\sc G.~Meunier}, {\em Homogenization for periodical electromagnetic structure:
  Which formulation?}, IEEE Transactions on Magnetics, 42 (2010), pp.~763--766.

\bibitem{monk2003finite}
{\sc P.~Monk}, {\em Finite element methods for Maxwell's equations}, Oxford
  University Press, 2003.

\bibitem{murat-gconv-77}
{\sc F.~Murat and L.~Tartar}, {\em {H}-convergence}.
\newblock S\'{e}minaire d'analyse fonctionnelle et num\'{e}rique de
  l'universit\'{e} d'Alger, 1977.

\bibitem{nguetseng-tsh-89}
{\sc G.~Nguetseng}, {\em A general convergence result for a functional related
  to the theory of homogenization}, S.I.A.M. Journal of Mathematical Analysis,
  20 (1989), pp.~608--623.

\bibitem{niyonzima-chm-16}
{\sc I.~Niyonzima, C.~Geuzaine, and S.~Sch{\"o}ps}, {\em Waveform relaxation
  for the computational homogenization of multiscale magnetoquasistatic
  problems}, Journal of Computational Physics, 327 (2016), pp.~416--433.

\bibitem{niyonzima-chm-12}
{\sc I.~Niyonzima, R.~V. Sabariego, P.~Dular, and C.~Geuzaine}, {\em Finite
  element computational homogenization of nonlinear multiscale materials in
  magnetostatics}, IEEE Transactions on Magnetics, 48 (2012), pp.~587--590.

\bibitem{niyonzima-chm-14}
{\sc I.~Niyonzima, R.~V. Sabariego, P.~Dular, and C.~Geuzaine}, {\em Nonlinear
  computational homogenization method for the evaluation of eddy currents in
  soft magnetic composites}, IEEE Transactions on Magnetics, 50 (2014),
  pp.~7001304, 1--4.

\bibitem{niyonzima-chm-13}
{\sc I.~Niyonzima, R.~V. Sabariego, P.~Dular, F.~Henrotte, and C.~Geuzaine},
  {\em Computational homogenization for laminated ferromagnetic cores in
  magnetodynamics}, IEEE Transactions on Magnetics, 49 (2013), pp.~2049--2052.

\bibitem{niyonzima-chm-13-1}
{\sc I.~Niyonzima, R.~V. Sabariego, P.~Dular, F.~Henrotte, and C.~Geuzaine},
  {\em A computational homogenization method for the evaluation of eddy current
  in nonlinear soft magnetic composites}, in Proceeding of the $9^\text{th}$
  International Symposium on Electric and Magnetic Fields, (EMF2013), Bruges,
  Belgium, April 2013.

\bibitem{pankov-gconv-97}
{\sc A.~Pankov}, {\em G-convergence and homogenization of nonlinear partial
  differential operators}, Kluwer academic publishers, 1997.

\bibitem{ren-hybrid-90}
{\sc Z.~Ren, F.~Bouillault, A.~Razek, A.~Bossavit, and J.-C. V{\'e}rit{\'e}},
  {\em A new hybrid model using electric field formulation for {3D} eddy
  current problems}, IEEE Transactions on Magnetics, 26 (1990), pp.~470--473.

\bibitem{rockafellar-convexanalysis-69}
{\sc R.~T. Rockafellar}, {\em Convex analysis}, Princeton Univ. Press,
  Princeton, NJ, 1969.

\bibitem{rodriguez-mqs-10}
{\sc A.~A. Rodr{\'\i}guez and A.~Valli}, {\em Eddy Current Approximation of
  Maxwell Equations: Theory, Algorithms and Applications}, vol.~4, Springer
  Science \& Business Media, 2010.

\bibitem{12.cefc.sabariego.homogLamHyst}
{\sc R.~V. Sabariego, I.~Niyonzima, C.~Geuzaine, and J.~Gyselinck}, {\em
  {Time-domain finite-element modelling of laminated iron cores -- Large skin
  effect homogenization considering the {Jiles-Atherton} hysteresis model}}, in
  Proceedings of the 15th Biennial IEEE Conference on Electromagnetic Field
  Computation (CEFC2012), Oita, Japan, November 11--14, 2012.

\bibitem{sanchez-homogenization-87}
{\sc E.~Sanchez-Palencia and A.~Zaoui}, {\em { Homogenization techniques for
  composite media}}, in Homogenization Techniques for Composite Media,
  vol.~272, 1987.

\bibitem{Sihvola-mixing-99}
{\sc A.~Sihvola}, {\em Electromagnetic mixing formulas and applications}, IEEE
  Electromagnetic Waves Series, 47), 1999.

\bibitem{tartar-homog-09}
{\sc L.~Tartar}, {\em The general theory of homogenization a personalized
  introduction}, Springer Berlin Heidelberg, 2009.

\bibitem{visintin-tsh-06-b}
{\sc A.~Visintin}, {\em Homogenization of doubly-nonlinear equations}, Rend.
  Lincei Mat. Appl., 17 (2006), pp.~211--222.

\bibitem{visintin-tsh-07-b}
{\sc A.~Visintin}, {\em Two-scale convergence of some integral functionals},
  Calc. Var., 29 (2007), pp.~239--265.

\bibitem{visintin-tsh-08}
{\sc A.~Visintin}, {\em Electromagnetic processes in doubly-nonlinear
  composites}, Communications in Partial Differential Equations, 33 (2008),
  pp.~804--841.

\bibitem{visintin-tsh-11}
{\sc A.~Visintin}, {\em Homogenization of a parabolic model of ferromagnetism},
  Journal of Differential Equations, 250 (2011), pp.~1521--1552.

\end{thebibliography}


\clearpage
\end{document}